\documentclass[11pt]{article}

\usepackage[margin=1.1in,a4paper]{geometry}
\usepackage[utf8]{inputenc}
\usepackage[T1]{fontenc}    % use 8-bit T1 fonts
\usepackage[colorlinks=true,citecolor=blue,linkcolor=blue, urlcolor=black]{hyperref}      % hyperlinks
\usepackage{url}            % simple URL typesetting
\usepackage{booktabs}       % professional-quality tables
\usepackage{amsfonts}       % blackboard math symbols
\usepackage{amsthm}
\usepackage{amsmath}
\usepackage{amssymb}
\usepackage{mdframed}
\usepackage{nicefrac}       % compact symbols for 1/2, etc.
\usepackage{microtype}      % microtypography
\usepackage{mathtools}
\usepackage{xcolor}         % colors
\usepackage{todonotes}
\usepackage{tikz}
\usepackage{tikz-cd}
\usepackage{framed}
\usepackage{enumitem}
\usepackage{soul}
\usepackage{cancel}
\newcommand{\subscript}[2]{$#1 _ #2$}

\mathtoolsset{showonlyrefs}

\newcommand{\sca}[2]{\langle #1, #2\rangle}
\newcommand{\nr}[1]{\left\Vert #1\right\Vert}
\newcommand{\abs}[1]{\left\vert #1\right\vert}
\newcommand{\Nsp}{\mathbb{N}}
\newcommand{\Rsp}{\mathbb{R}}

\newcommand{\Class}{\mathcal{C}}

\newcommand{\X}{\mathcal{X}}
\newcommand{\Y}{\mathcal{Y}}

\newcommand{\dd}{\mathrm{d}}

\newcommand{\Lip}{\mathrm{Lip}}

\newcommand{\eps}{\varepsilon}

\newcommand{\KL}{\mathrm{KL}}
\newcommand{\Prob}{\mathcal{P}}

\newcommand{\Var}{\mathbb{V}\mathrm{ar}}

\newcommand{\Esp}{\mathbb{E}}
\newcommand{\Kant}{\mathcal{K}}

\newcommand{\spt}{\mathrm{spt}}

\newcommand{\kl}[2]{\ensuremath{\mathrm{KL}\left(#1 \vert #2\right)}}

\newmdtheoremenv{theo}{Theorem}
\newmdtheoremenv{coro}{Corollary}

\newtheorem{theorem}{Theorem}[section]

\newtheorem{corollary}[theorem]{Corollary}
\newtheorem{lemma}[theorem]{Lemma}

\newtheorem{proposition}[theorem]{Proposition}

\theoremstyle{definition}
\newtheorem{definition}{Definition}[section]

\newtheorem{remark}{Remark}[section]

\title{Sharper Exponential Convergence Rates for Sinkhorn's Algorithm in Continuous Settings}

\author{
Lénaïc Chizat\thanks{Ecole Polytechnique Fédérale de Lausanne (EPFL), Institute of Mathematics, 1015 Lausanne, Switzerland. Authors ordered alphabetically.}
\and Alex Delalande\footnotemark[1] \textsuperscript{,}\thanks{Corresponding author: alex.delalande@epfl.ch} 
\and Tomas Vaškevičius\footnotemark[1]}

\begin{document}

\maketitle
%\footnotetext[1]{}

\begin{abstract}
  We study the convergence rate of Sinkhorn's algorithm for solving entropy-regularized optimal transport problems when at least one of the probability measures, $\mu$, admits a density over $\mathbb{R}^d$. For a semi-concave cost function bounded by $c_{\infty}$ and a regularization parameter $\lambda > 0$, we obtain exponential convergence guarantees on the dual sub-optimality
  gap with contraction rates that are polynomial in $\lambda/c_{\infty}$. This represents an exponential improvement over
  the known contraction rate $1 - \Theta(\exp(-c_{\infty}/\lambda))$ achievable via
  Hilbert's projective metric.
  Specifically, we prove a contraction rate value of $1-\Theta(\lambda^2/c_\infty^2)$ when $\mu$ has a bounded log-density. In some cases, such as when $\mu$ is log-concave and the cost function is $c(x,y)=-\sca{x}{y}$, this rate improves to $1-\Theta(\lambda/c_\infty)$. The latter rate matches the one that we derive for the transport between isotropic Gaussian measures, indicating tightness in the dependency in $\lambda/c_\infty$. Our results are fully non-asymptotic and explicit in all the parameters of the problem.
\end{abstract}

{
}

\section{Introduction}
\label{sec:introduction}

We study the numerical resolution of the
entropy-regularized optimal transport problem
\begin{equation}
  \label{eq:entropic-ot-primal}
  \inf_{\gamma \in \Gamma(\mu, \nu)}
  \int
  c(x,y)\dd\gamma(x,y)
  + \lambda \kl{\gamma}{\mu \otimes \nu}.
\end{equation}
Here, $\mu$ and $\nu$ are probability measures 
supported over subsets
$\X$ and $\Y$ of $\mathbb{R}^{d}$, 
$\Gamma(\mu, \nu)$ denotes the set of couplings between $\mu$ and $\nu$,
$\lambda$ is a positive regularization parameter,
$\mathrm{KL}$ is the Kullback-Leibler divergence
and $c(x,y)$ is the cost of assigning a unit of mass from $x$ to $y$.
Throughout this paper, we assume that the cost $c$ has bounded 
oscillation 
on $\mathcal{X} \times \mathcal{Y}$ and we denote
\begin{equation}
  \label{eq:boundedness-condition}
  c_{\infty} = \nr{c}_{\mathrm{osc}} = \sup_{x \in \mathcal{X}, y \in \mathcal{Y}} c(x,y) - \inf_{x \in \mathcal{X}, y \in \mathcal{Y}} c(x,y) < \infty.
\end{equation}
Under the conditions listed above, the unique solution to
problem \eqref{eq:entropic-ot-primal} can be obtained by solving its dual problem
\begin{equation}
  \label{eq:dual-objective}
  \sup_{\phi \in L^{  \infty  }(\mu), \psi \in L^{  \infty  }(\nu)} \int \phi \dd\mu + \int \psi \dd\nu
  + \lambda \left(
    1 -
    \int \int \exp\left(\frac{\phi \oplus \psi - c}{\lambda}\right)
    \dd\mu \dd\nu 
  \right),
\end{equation}
where $(\phi \oplus \psi) (x,y) = \phi(x) + \psi(y)$. The classical approach to solve
\eqref{eq:dual-objective} is via \emph{Sinkhorn's algorithm}. Denoting 
$F(\phi, \psi)$ the dual functional being maximized in \eqref{eq:dual-objective}, 
Sinkhorn's algorithm starts at an arbitrary point
$\psi_{0} \in L^{  \infty  }(\nu)$ and 
repeatedly performs exact maximization of the objectives
$\phi \mapsto F(\phi, \psi)$
and $\psi \mapsto F(\phi, \psi)$. This leads to the sequence of iterates
\begin{equation}
  \label{eq:sinkhorn-iterates-informal}
  \begin{cases}
    \phi_{t + \frac{1}{2}} \in \arg \max_{\phi} F(\phi, \psi_{t}),
    \\
    \psi_{t+1} \in \arg \max_{\psi} F(\phi_{t+\frac{1}{2}}, \psi),
  \end{cases}
\end{equation}
which can be defined in a unique way as detailed in
Section~\ref{sec:background}. For Sinkhorn's iterates $(\phi_{t+1/2}, \psi_{t})_{t \geq 0}$
and any integer $t \geq 0$, we denote the dual sub-optimality gap
\begin{equation}
  \delta_{t} = \sup_{\phi \in L^{  \infty  }(\mu),
  \psi \in L^{  \infty  }(\nu)}
  F(\phi, \psi) - F(\phi_{t+1/2}, \psi_{t}).
\end{equation}
It is well-known that $\delta_{t}$ converges to zero under general
conditions. However, despite the extensive history of Sinkhorn's algorithm and its
widespread adoption in modern applications, several open questions remain regarding
its convergence speed.

Existing guarantees generally fall into two categories. The first category
exhibits exponential convergence, albeit with an exponential
dependence on $c_{\infty}/\lambda$ in the rate.
%, namely$\delta_{t} = O(c_{\infty}(1-e^{-c_{\infty}/\lambda})^{t})$.
Although these guarantees demonstrate the exponential nature of 
the convergence of Sinkhorn's algorithm, they often fail to provide bounds that are relevant in
practical applications where $c_{\infty}/\lambda$ is large.
This limitation is somewhat mitigated by the second category
of guarantees that yield bounds of the form
$\delta_{t} = O(c_{\infty}^{2}/(\lambda t))$.
While these bounds significantly improve the dependence on $c_{\infty}/\lambda$, they exhibit a slower, polynomial dependence on the number
of iterations $t$. We discuss both types of guarantees in greater depth in
Section~\ref{sec:related-work}.

%\tomas{double-check if the dependence on $c_{\infty}$ in the above rates is correct}

One feature of the 
two types of convergence guarantees discussed above is their
generality, as they require no assumptions on the marginal measures $\mu$ and $\nu$.
However, one of the chief advantages of the optimal transport problem
\eqref{eq:entropic-ot-primal} lies in its capacity to capture specific geometric properties of the problem at hand, which underscores its importance in a broad array of application domains. This leads to the main question addressed in our paper: can we improve upon the convergence guarantees described above under additional geometric assumptions on the cost
function $c$ and the marginal measures $\mu$ and $\nu$? 

Our main results demonstrate that the presence of geometric structure within the
problem can be leveraged to obtain exponential convergence rates with a
well-behaved dependence on $\lambda/c_{\infty}$.
In particular, we prove convergence bounds of the form
$(1 - \Theta(\lambda/c_{\infty}))^{t}$ and $(1-\Theta(\lambda^2/c_{\infty}^{2}))^{t}$,
representing a best-of-both-worlds scenario in terms of the bounds
discussed above.
For the sake of simplicity and for streamlining the comparison with related
work in Section~\ref{sec:related-work}, we begin by stating a special case of
our main results for the linear cost $c(x,y) = -\langle x, y \rangle$, which is
equivalent to the quadratic
cost $c(x,y)=\nr{x-y}^2$. 
The general statements and proofs of our results, which cover for instance 
any $\mathcal{C}^2$ cost over a bounded domain, are deferred to 
Section~\ref{sec:main-results}. 
We split our results into two cases, depending on the assumptions imposed on
one of the marginal measures, keeping the other marginal assumptions-free.
First, we handle the setting where one of the marginals is log-concave.
\begin{theorem}
  \label{thm:quadratic-cost-log-concave}
  Let $\X, \Y$ be compact with $\X$ convex. Let $c(x,y) = -\langle x, y
  \rangle$.
  If $\mu$ admits a log-concave density on $\mathcal{X}$, then 
  for any $t \geq 0$ we have
  $
    \delta_{t}
    \leq
    \delta_{0}
    \left(1 - \alpha^{-1}\right)^{t}$ 
    where $\alpha = 176( 1 + \frac{c_{\infty}}{\lambda})$.
\end{theorem}
Second, we analyze convergence under the assumption that one of the marginal
measures has a bounded log-density.
\begin{theorem}
  \label{thm:quadratic-cost-bounded-density}
  Let $\X, \Y$ be compact with $\X$ convex. Let $c(x,y) = -\langle x, y
  \rangle$. Suppose that $\mu$ admits a density $f$ satisfying $0 < m \leq f \leq M < \infty$ for
  some positive constants $m$ and $M$. Then:
  \begin{enumerate}
    \item For any $t \geq 0$ we have
      $\delta_{t} \leq \delta_{0}(1 - \alpha^{-1})^t$,
      where $\alpha = 176(1 + \frac{M}{m}\frac{c_{\infty}}{\lambda} +
      \frac{c_{\infty}^2}{\lambda^2})$.
    \item For any large enough $t$ we have
      $\delta_{t+1} \leq \delta_{t}(1 - \alpha^{-1})$,
      where $\alpha = 176(1 + \frac{M}{m} \frac{c_{\infty}}{\lambda})$.
  \end{enumerate}
\end{theorem}%
Finally, we prove a matching lower bound for the contraction rates proved in Theorem~\ref{thm:quadratic-cost-log-concave}
and the second part of Theorem~\ref{thm:quadratic-cost-bounded-density}. Our lower bound holds for one-dimensional zero-mean Gaussian measures.
While the Gaussian measures do not satisfy the compactness assumptions considered in Theorems~\ref{thm:quadratic-cost-log-concave} and \ref{thm:quadratic-cost-bounded-density}, our lower bound provides strong supporting evidence for the optimality of the $1 - \Theta(\lambda/c_{\infty})$ contraction rate.
The proof is deferred to Section~\ref{sec:gaussians}, where the limiting contraction rate of the sequence $(\delta_t)_{t \geq 0}$ is also given in the asymptotic regime $\lambda \to 0$.
\begin{theorem}
    \label{thm:lower-bound-gaussian-less-formal}
    On $\Rsp$, let $\mu = \mathcal{N}(0, 1)$, $\nu = \mathcal{N}(0, \sigma^2)$ for $\sigma > 0$, $c(x,y) = -xy$, and $\psi_0 = 0$.
    If $\lambda \leq \sigma/5$, then for any $t \geq 0$ we have
    $\delta_{t} \geq \frac{\sigma}{20} (1 - \frac{5\lambda}{\sigma})^{t}$.
\end{theorem}

\subsection{Comparison with related work}
\label{sec:related-work}

Iterations~\eqref{eq:sinkhorn-iterates-informal}, which we refer to as Sinkhorn's algorithm \cite{sinkhorn1964algorithm}, were rediscovered several times and perhaps first appeared in~\cite{yule1912methods}. Consequently, these iterations and close variants thereof are known under various names in various communities, such as the Sinkhorn-Knopp algorithm, Iterative Proportional Fitting Procedure (IPFP), RAS algorithm, matrix scaling algorithm, Bregman algernative projection or Fortet's iterations; see the surveys~\cite{idel2016review, peyre2019computational} for some historical remarks. Most quantitative convergence analyses can be classified into two categories: (i)~\emph{exponential} convergence rates with \emph{non-robust} constants and (ii)~\emph{polynomial} convergence rates with \emph{robust} constants.

Convergence guarantees of the first kind can be obtained by noticing that Sinkhorn's iteration is a contraction for Hilbert's projective metric~\cite{franklin1989scaling, borwein1994dad} where the contraction rate can be shown to be at most $1-\Theta(e^{-c_\infty/\lambda})$; see~\cite{chen2016hilbertmetric, deligiannidis2024hilbertmetric} for extensions to continuous settings and \cite{eckstein2024hilbertsprojectivemetricfunctions} for non-compact settings, as well as~\cite[Remark~4.12]{peyre2019computational} for a bibliographical discussion. Note that this type of contraction rate can also be obtained without using Hilbert's projective metric, for example by using a classical convex optimization approach as in \cite{marino2020schrodinger, carlier2022multimarginalsinkhorn}, allowing to extend the analysis to multi-marginal settings; or by relying on coupling methods and stochastic control following \cite{greco2023coupling, conforti2023quantitative}, allowing extensions to non-bounded configurations. As put forth in \cite[Remark~4.15]{peyre2019computational}, the dependence on $c_\infty/\lambda$ of the $1-\Theta(e^{-c_\infty/\lambda})$ contraction rates numerically seems tight in the most difficult cases, e.g. when the cost $c$ is random. In more favorable cases however, such as when there exists an optimal transport map solving problem \eqref{eq:entropic-ot-primal} with $\lambda=0$, many have observed experimentally that the dependency in $c_\infty/\lambda$ is generally much better in practice.
%at least in ``geometric'' settings~\cite[Remark.~4.15]{peyre2019computational}. Some recent works have extended such guarantees to the non-compact setting~\cite{conforti2023quantitative}. 

Convergence guarantees of the second kind can be traced back to~\cite{kalantari2008complexity}, and have more recently regained interest~\cite{altschuler2017near, chakrabarty2021sinkhornsublinearrate, dvurechensky2018computational}. In their tightest form, due to~\cite{dvurechensky2018computational}, one has the guarantee $\delta_t \leq 2c_\infty^2/(\lambda t)$ for $t\geq 1$ (see~\cite[Proposition~10]{chizat2020faster}). This result relies on an explicit expression (see Eq.~\eqref{eq:sinkhorn-iterates-improvement}) for the one-step dual improvement $\delta_{t}-\delta_{t+1}$, which is also central in our approach. This proof scheme can also be framed as a mirror descent analysis in well-chosen geometries~\cite{leger2021gradient,aubin2022mirror} and has been recently extended to non-compact settings~\cite{ghosal2022convergence} or~\cite[Proposition~4]{de2021diffusion}.

%\begin{itemize}
%\item For (i): Hilbert's metric technique~\cite{franklin1989scaling}, see~\cite[Remark.~4.2]{peyre2019computational} for bibliographical discussion.  Extension to the non-compact case with different techniques~\cite{conforti2023quantitative}. Typically smth like $(1-e^{\Vert c\Vert_\infty/\lambda})^t$.
%\item For (ii): proof technique introduced by~\cite{altschuler2017near}, then improved by~\cite{dvurechensky2018computational} (obtains $\propto \Vert c\Vert_\infty/(\lambda\cdot t)$); link with mirror descent proof~\cite{leger2021gradient,aubin2022mirror}; extended to non-compact case in~\cite{ghosal2022convergence} or~\cite[Prop.~4]{de2021diffusion}.
%\end{itemize}
Our results in this paper take the best of both worlds since, under additional regularity assumptions on the cost and marginals, we prove \emph{exponential} convergence rates with \emph{robust} constants. Some prior works~\cite{ berman2020sinkhorn,deb2023wasserstein} have studied, for the cost $c(x,y)=-\langle x,y \rangle$, the limiting dynamics of the time-rescaled sequence of iterates $(\psi_{t/\lambda})_{t}$ when $\lambda\to 0$. Under appropriate assumptions, this leads to a dynamics that converges at an (unspecified) exponential rate to the unregularized dual potentials. Informally, this might be seen as a hint towards a convergence rate in $1-\Theta(\lambda/c_\infty)$ in certain settings, which is indeed one of our results.  

Finally, a result related to ours was announced by S. Di Marino in the 2022 workshop \emph{Optimal Transportation and Applications} in Pisa. Specifically, he announced a dual convergence rate of the form $1-\Theta(\lambda^2/c_\infty^2)$ for bounded and Lipschitz cost functions when both measures $\mu$ and $\nu$ satisfy a Poincaré inequality\footnote{Personal communication.}, which differs from the assumptions considered in our work. S. Di Marino's proof extends to the multimarginal setting, but it requires assumptions on all the marginals, while we only make assumptions on one of the marginals.

%The works which are more directly related are:
%\begin{itemize}
%\item the small regularization asymptotic of Sinkhorn, which implies a rate in $(1-C\lambda)$ (but $C>0)$ is not explicit I think)  
%\item Simone's ghost work
%\end{itemize}
%In contrast to the above mentionned works, our result is non-asymptotic in $\lambda$ and in $t$, with fully explicit convergence bounds.

%\begin{equation}
%  \kl{\gamma}{\mu \otimes \nu} = \int
%  \log \frac{\d\gamma}{\dd(\mu \otimes \nu)}
%  \dd\gamma(x,y).
%\end{equation}

\subsection{Paper outline}
\label{sec:paper-outline}
The rest of the paper is organized as follows. In Section~\ref{sec:background} we define the objects of interest and recall a few useful formulas. The full statements of our main results are in Section~\ref{sec:main-results}, which also contains an outline of the proofs in Section~\ref{sec:exponential-convergence-thm-proof}. The core of the proofs can be found in Sections~\ref{sec:proof-of-one-step-improvement-proposition} and \ref{sec:proof-of-variance-subopt-bound}. The convergence rate of Sinkhorn's algorithm for isotropic Gaussian is studied in Section~\ref{sec:gaussians}. Section~\ref{sec:numerical-experiments} gathers numerical experiments indicating that the core feature of our theoretical findings also holds in practical settings where both marginals are discrete probability distributions. Finally we gather in the Appendix the proof of intermediate technical results.

\section{Background and notation}
\label{sec:background}

This section sets up the notation and provides background
on the Sinkhorn algorithm and the (semi-)dual entropic optimal transport
objective.

\subsection{General conventions}
For a subset $\mathcal{Z}$ of some Euclidean space we let
$\mathcal{P}(\mathcal{Z})$ be the set of all Borel probability measures
supported on $\mathcal{Z}$. For any $\rho \in \Prob(\mathcal{Z})$ and any $f \in L^2(\rho)$, the expectation and variance of $f$ against the 
measure $\rho$ are respectively denoted 
\begin{equation*}
    \Esp_\rho f= \Esp_{x \sim \rho} f(x) = \int f \dd \rho \quad \text{and} \quad \Var_\rho(f) = \Var_{x \sim \rho} (f(x)) = \Esp_\rho (f - \Esp_\rho f)^2.
\end{equation*}
For any $\rho, \rho' \in \mathcal{P}(\mathcal{Z})$,
$\kl{\rho}{\rho'}$ is the Kullback-Leibler divergence (or relative
entropy) between the measures $\rho$ and $\rho'$, defined as
\begin{equation}
  \kl{\rho}{\rho'}
  = \begin{cases}
    \int \log\frac{\dd \rho}{\dd \rho'} \dd \rho & \text{ if } \rho \ll \rho', \\
    +\infty & \text{ otherwise.}
  \end{cases}
\end{equation}
Throughout the paper, the norm $\|\cdot\|$ is the Euclidean norm,
$\mathcal{X} \subseteq \mathbb{R}^{d}$ is the support of the marginal measure
$\mu \in \mathcal{P}(\mathcal{X})$ and $\mathcal{Y} \subseteq \mathbb{R}^{d}$
is the support of the marginal measure $\nu \in \mathcal{P}(\mathcal{Y})$.
For functions $\phi \in L^{1}(\mu)$ and $\psi \in L^{1}(\nu)$ we use
the notation
\begin{equation}
  \sca{\phi}{\mu} = \int \phi \dd \mu
  \quad\text{and}\quad
  \sca{\psi}{\nu} = \int \psi \dd \nu.
\end{equation}
We denote the set of real-valued continuous
functions defined on $\mathcal{Z}$ by $\mathcal{C}(\mathcal{Z})$. For a function $f : \mathcal{Z} \to \mathbb{R}$, we denote
its oscillation norm by $\|f\|_{\mathrm{osc}}
= \sup_{x \in \mathcal{Z}} f(z) - \inf_{x' \in \mathcal{Z}} f(z)$ and 
its Lipschitz semi-norm by $\Lip(f) = \sup_{(x,y) \in \mathcal{Z}\times\mathcal{Z}} (f(y) - f(x))/\nr{y-x}$. For the cost function $c$,
we use the notation $c_{\infty} = \|c\|_{\mathrm{osc}}$.
Finally, we recall the notion of semi-concave and semi-convex functions that appears throughout our results:
\begin{definition}[Semi-concave/convex functions]
    On a convex subset $\mathcal{Z}$ of some Euclidean space, a function $f : \mathcal{Z} \to \Rsp$ is said to be $\xi$-semi-concave for some $\xi \in \Rsp_+$ if the function $x \mapsto f(x) - \frac{\xi}{2}\nr{x}^2$ is a concave function. Similarly, $f$ is said to be $\zeta$-semi-convex for some $\zeta \in \Rsp_+$ if $x \mapsto f(x) + \frac{\zeta}{2}\nr{x}^2$ is a convex function.
\end{definition}

\subsection{\texorpdfstring{$(c,\lambda)$}{(c,lambda)}-transforms and the
Sinkhorn algorithm}

For a cost function $c$ and a regularization parameter $\lambda > 0$,
the $(c,\lambda)$-transforms $\phi^{c,\lambda}$,$\psi^{c,\lambda}$
of functions $\phi \in L^{  \infty  }(\mu), \psi \in L^{  \infty  }(\nu)$
are respectively defined as
\begin{align}
  \label{eq:c-lambda-transform}
  \phi^{c, \lambda}(y)
  &= -\lambda \log \int \exp \left(\frac{\phi(x) -
  c(x,y)}{\lambda}\right)\dd\mu(x),  \\ %\quad \phi^{c,\lambda} \in  L^{  \infty  }(\nu),
  \psi^{c, \lambda}(x)
  &= -\lambda \log \int \exp \left(\frac{\psi(y) -
  c(x,y)}{\lambda}\right)\dd\nu(y). %,\quad \psi^{c,\lambda} \in L^{  \infty  }(\mu)
\end{align}
Note that even thought $\phi$ and $\psi$ are only defined $\mu$- and $\nu$-almost-everywhere, 
their transforms $\phi^{c, \lambda}$ and $\psi^{c, \lambda}$ are defined everywhere on $\Rsp^d$.
A direct computation shows that if $c$ is a function of bounded oscillation or is Lipschitz continuous, 
so are the $(c, \lambda)$-transforms:
\begin{lemma}
\label{lemma:properties-c-transforms}
    If $\nr{c}_{\mathrm{osc}} = c_\infty < \infty$ (resp. $\Lip(c)<\infty)$, then for any $\phi \in L^{  \infty  }(\mu)$ and $\psi \in L^{  \infty  }(\nu)$,
    \begin{align}
      \label{eq:c-lambda-oscilation}
      \|\phi^{c,\lambda}\|_{\mathrm{osc}}
      \leq c_{\infty}
      \quad &\text{and}\quad
      \|\psi^{c,\lambda}\|_{\mathrm{osc}} \leq c_{\infty}, \\
      (\text{resp.} \quad 
      \Lip(\phi^{c,\lambda}) \leq \Lip(c)
      \quad &\text{and}\quad
      \Lip(\psi^{c,\lambda}) \leq \Lip(c)).
      %\quad\text{for any}\quad
      %\phi \in L^{  \infty  }(\mu), \psi \in
      %L^{  \infty  }(\nu).
    \end{align}
\end{lemma}
The $(c,\lambda)$-transform can be seen as an entropic equivalent of the
$c$-transform in the unregularized optimal transport theory.
Indeed, it can be shown by computing the first variation of
the dual objective \eqref{eq:dual-objective} that
for any $\phi \in L^{  \infty  }(\mu)$ and $\psi \in
L^{  \infty  }(\nu)$, we have
\begin{align}
  F(\psi^{c,\lambda}, \psi) = \max_{\phi \in L^{  \infty  }(\mu)} F(\phi,
  \psi)
  \quad\text{and}\quad
  F(\phi, \phi^{c,\lambda})
  =
  \max_{\psi \in L^{  \infty  }(\nu)} F(\phi, \psi).
\end{align}
For convenience, we also introduce the following shorthand notation for
applying the $(c,\lambda)$-transformation twice:
\begin{equation}
  \label{eq:c-lambda-transform-twice}
  \phi^{\overline{c,\lambda}} = (\phi^{c,\lambda})^{c, \lambda}
  \quad\text{and}\quad
  \psi^{\overline{c,\lambda}} = (\psi^{c,\lambda})^{c,\lambda}.
\end{equation}
A direct computation shows that
a maximizing pair $(\phi^{*}, \psi^{*})$ for the dual objective
$F : L^{  \infty  }(\mu) \times L^{  \infty  }(\nu) \to \mathbb{R}$
can be obtained by solving the Schr\"{o}dinger system:
\begin{equation}
  \label{eq:schroedinger-system}
  \begin{cases}
    \phi(x) = \psi^{c,\lambda}(x)\quad \text{for every } x \in \mathcal{X},\\
    \psi(y) = \phi^{c,\lambda}(y)\quad\text{for every }y\in\mathcal{Y}.
  \end{cases}
\end{equation}
In particular, a solution $(\phi, \psi)$ to the Schr\"{o}dinger system verifies
$\phi = \psi^{c,\lambda}$ and $\psi = \psi^{\overline{c,\lambda}}$.
Under our assumptions, such a maximizing pair always exists and it is unique up
to translations by a constant.
Hence, we may now more formally define the Sinkhorn updates
initialized at an arbitrary $\psi_{0} \in L^{  \infty  }(\nu)$ by
\begin{align}
  \begin{split}
    \label{eq:sinkhorn-iterates}
    \phi_{t+1/2}
    = \psi^{c,\lambda}_{t}
    \quad\text{and}\quad
    \psi_{t+1} = \phi_{t+1/2}^{c,\lambda} = \psi_{t}^{\overline{c,\lambda}}
    \quad\text{for }t \geq 0.
  \end{split}
\end{align}
Because for every $t \geq 0$ the iterate $(\phi_{t+1/2}, \psi_{t}) =
(\psi_{t}^{c,\lambda}, \psi_{t})$ depends only on $\psi_{t}$,
it will be convenient to rewrite the dual objective
in terms of the function $\psi \in L^{  \infty  }(\nu)$.
This leads to the notion of semi-dual objective described in the next section.

\subsection{Semi-dual objective and primal-dual relations}
The semi-dual objective $E : L^{  \infty  }(\nu) \to \mathbb{R}$ is
defined by
\begin{equation}
  \label{eq:semi-dual}
  E(\psi) =
  \max_{\phi \in L^{  \infty  }(\mu)}
  F(\phi,\psi)
  = F(\psi^{c,\lambda}, \psi).
\end{equation}
The definition of the $(c,\lambda)$-transforms implies that for any
$\psi \in L^{  \infty  }(\nu)$
\begin{equation}
  \int\int \exp\left(\frac{\psi^{c,\lambda} \oplus \psi - c}{\lambda}\right) \dd\mu \dd\nu
  = 1.
\end{equation}
Hence, for any $\psi \in L^{  \infty  }(\nu)$ the semi-dual objective $E$
can be written as
\begin{equation}
  E(\psi) = \int \psi^{c,\lambda} \dd \mu + \int \psi \dd \nu.
\end{equation}
Any $\psi \in L^{  \infty  }(\nu)$
induces a probability measure $\gamma[\psi] \in \mathcal{P}(\mathcal{X} \times
\mathcal{Y})$ defined as
\begin{equation}
  \label{eq:dual-to-primal}
  \dd\gamma[\psi](x, y) = \exp\left(\frac{\psi^{c,\lambda}(x) + \psi(y) -
  c(x,y)}{\lambda}\right)\dd\mu(x)\dd\nu(y).
\end{equation}
The measure $\gamma[\psi]$ is, in general, not primal feasible as it does not
belong to the set of couplings $\Gamma(\mu, \nu)$. In fact, we have that 
$\gamma[\psi]$ belongs to $\Gamma(\mu,\nu)$ if and only if $\psi$ is equal to a maximizer $\psi^*$ of the semi-dual objective $E$, in which case 
the measure $\gamma[\psi] = \gamma[\psi^{*}]$ is equal to the unique minimizer
of the primal problem \eqref{eq:entropic-ot-primal}. Moreover,
for any $\psi \in L^{  \infty  }(\nu)$ the following relation connects the
dual sub-optimality gap to a form of a primal guarantee:
\begin{align}
  E(\psi^{*}) - E(\psi)
  = \lambda \kl{\gamma[\psi^{*}]}{\gamma[\psi]}.
\end{align}

\subsection{Basic properties of Sinkhorn's iterates}
We now list some basic properties satisfied by Sinkhorn's iterates
$(\phi_{t+1/2}, \psi_{t})_{t \geq 0}$. They are well-known in the
literature, and we will use them without explicit references to this
section.
For any $\psi \in L^{  \infty  }(\nu)$, define the marginal measures
$\mu[\psi] \in \mathcal{P}(\mathcal{X})$ and $\nu[\psi] \in
\mathcal{P}(\mathcal{Y})$
of $\gamma[\psi] \in \mathcal{P}(\mathcal{X} \times \mathcal{Y})$ by
\begin{align}
  \label{eq:marginal-measures-definition}
  \mu[\psi](A) &= \int_{A \times \mathcal{Y}} d\gamma[\psi](x,y)
  \quad\text{for any Borel set }A \subseteq \mathcal{X},
  \\
  \nu[\psi](B) &= \int_{\mathcal{X} \times B} d\gamma[\psi](x,y)
  \quad\text{for any Borel set }B \subseteq \mathcal{Y}.
\end{align}
The definition of $(c,\lambda)$-transforms ensures that
$\mu[\psi] = \mu$.
However, $\nu[\psi] \neq \nu$ unless $\psi$ is a maximizer of the semi-dual
objective $E$. An explicit calculation shows that
for any $\psi \in L^{  \infty  }(\nu)$,
\begin{equation}
  \psi^{\overline{c,\lambda}} - \psi
  = \lambda \log \frac{\dd \nu}{\dd \nu[\psi]},
\end{equation}
so that 
\begin{equation}
\label{eq:KL-nu-nu-psi}
    \kl{\nu}{\nu[\psi]} = \frac{1}{\lambda}\sca{\psi^{\overline{c, \lambda}} - \psi}{\nu}.
\end{equation}
In particular, because $\psi_{t+1} = \psi_{t}^{\overline{c,\lambda}}$,
for any $t \geq 0$ we have
\begin{equation}
  \label{eq:sinkhorn-iterates-improvement}
  E(\psi_{t+1}) - E(\psi_{t})
  \geq
  F(\phi_{t+1/2}, \psi_{t+1}) - F(\phi_{t+1/2}, \psi_{t})
  = \lambda \kl{\nu}{\nu[\psi_{t}]}.
\end{equation}
Finally, for any $t \geq 1$ and any maximizer $\psi^{*}$ of the semi-dual objective $E$,
it holds following from \eqref{eq:c-lambda-oscilation} that 
$\|\psi_{t} - \psi^{*}\|_{\mathrm{osc}} \leq 2c_{\infty}$.

\subsection{Entropic Kantorovich functional}
We end this background section with the introduction of 
what we refer to as the \emph{entropic Kantorovich functional}, denoted 
$\Kant : L^{  \infty  }(\nu) \to \Rsp$ and defined by
\begin{equation}
  \Kant(\psi) : \psi \mapsto \sca{\psi^{c, \lambda}}{\mu}.
\end{equation}
Notice that this functional 
corresponds to the non-linear part of the semi-dual 
functional $E(\psi) = \sca{\psi^{c,\lambda}}{\mu} + \sca{\psi}{\nu} = \Kant(\psi)+ \sca{\psi}{\nu}$. The directional derivatives of $\Kant$ can be made explicit. First, for any $\psi \in L^{  \infty  }(\nu)$ and any $x \in \mathcal{X}$, we define the probability measure $\nu_x[\psi]$ on $\Y$ as
  \begin{align}
    \label{eq:nu-x-dfn}
      \dd \nu_x[\psi](y) &= \frac{e^{\frac{\psi(y) - c(x,y)}{\lambda}}}{ \int e^{\frac{\psi(\tilde{y}) - c(x,\tilde{y})}{\lambda}} \dd \nu(\tilde{y}) } \dd \nu(y) = e^{\frac{\psi^{c,\lambda}(x) + \psi(y) - c(x,y)}{\lambda}} \dd \nu(y) = \frac{\dd \gamma[\psi](x,y)}{\dd \mu(x)}.
  \end{align}
  Notice that the family $(\nu_x[\psi])_{x \in \X}$ corresponds to the disintegration of $\nu[\psi]$ with respect to $\mu$, in the sense that for any Borel set $B \subset \Y$ it holds
  $$ \nu[\psi](B) = \int \nu_x[\psi](B) \dd \mu(x). $$
We will repeatedly make use of the following property, proven in Section~\ref{sec:proof-of-derivatives-lemma}. % satisfied by the entropic Kantorovich functional.
\begin{lemma}
  \label{lemma:K-concave-derivatives}
    The functional $\Kant$ is concave on $L^{  \infty  }(\nu)$, and it is strictly concave up to translations.
    Moreover, $\Kant$ admits the following first and second order directional 
    derivatives: for any functions $\psi, v \in L^{  \infty  }(\nu)$ and any $\varepsilon \in \mathbb{R}$,
    \begin{align*}
        \frac{\dd}{\dd \eps} \Kant(\psi + \eps v) = -\sca{v}{\nu[\psi + \eps v]}, \quad \text{and} \quad 
        \frac{\dd^2}{\dd \eps^2} \Kant(\psi + \eps v) = -\frac{1}{\lambda}
        \int \Var_{\nu_x[\psi + \eps v]} (v) \dd \mu(x).
    \end{align*}
\end{lemma}

\section{Main result: exponential convergence with robust contraction constants}
\label{sec:main-results}

%\subsection{Statement and discussion}

Our main result establishes exponential convergence
of Sinkhorn's iterates in the semi-dual objective under various general conditions on the cost and one of the marginals: %whenever the cost function $c$ is uniformly semi-concave in its first variable.

  \begin{theorem}
    \label{thm:exponential-convergence}
    Assume that
    \begin{enumerate}[label=(\textit{A})]
        \item \label{assump} The domain $\X$ is convex, there exists $\xi \in \Rsp_+$ such that for all $y \in \Y$, $x \mapsto c(x, y)$ is $\xi$-semi-concave, and
      $\nr{c}_{\mathrm{osc}} = c_{\infty} < \infty$.
    \end{enumerate}
    Then, for any integer $t \geq 0$, the Sinkhorn iterates $(\psi_{t})_{t\geq 0}$ defined in
    \eqref{eq:sinkhorn-iterates} satisfy
    \begin{equation} \label{eq:exponential-convergence}
      %E(\psi^{*}) - E(\psi_{t}) \leq (1 - \alpha^{-1})^{t}(E(\psi^{*}) - E(\psi_{0}))
      E(\psi^{*}) - E(\psi_{t+1}) \leq (1 - \alpha^{-1})(E(\psi^{*}) - E(\psi_{t}))
    \end{equation}
    provided either one of the following additional assumption holds:
    \begin{enumerate}[label=(\subscript{A}{{\arabic*}})]
      \item \label{assump1} The domain $\X$ is compact and included in $\{ x : \|x\| \leq R_{\mathcal{X}}\}$, the measure $\mu$ admits a density $f(x)$ such that
        %\begin{equation}
          $
          \frac{
            \sup_{x \in \mathcal{X}} f(x)}{
            \inf_{x' \in \mathcal{X}} f(x')
          } = \kappa < \infty,
          $
        %\end{equation}
        %$t \geq 0$ is arbitrary,
        and
        $$ \alpha = 176 \{ 1 +
        (c_\infty + \frac{\xi}{2} R_\X^2)\kappa \lambda^{-1}
        + c_\infty^2 \lambda^{-2} \}.
        $$
        %Furthermore, for any $\mu,\nu,c$ and $\psi_{0}$ there exists an integer
        %$T$ such that for any $t \geq T$ the contraction constant $\alpha$
        %improves to
        %$$
        %  \alpha =
        %  96
        %  +
        %  96e\kappa (c_{\infty} +
        %  \frac{\xi}{2}R_{\mathcal{X}}^{2})\lambda^{-1}.
        %$$
        \item \label{assump2} There exists a $\xi$-strongly convex function $V
          : \X \to \Rsp$ such that the density of $\mu$ reads $f(x) =
          e^{-V(x)}$, and
        $$
        \alpha = 176\{ 1 +
        c_{\infty}  \lambda^{-1}
        +
        c_{\infty}^{2}\lambda^{-2}
        \}.
        $$
        %\tomas{Alex, could you add a short remark below the theorem addressing the fact that $\xi$ does not appear in $\alpha$? (The same also applies to the statement below under \ref{assump3}).} 
        \item \label{assump3} There exists $\zeta \in \Rsp_+$ such that for all
          $y \in \Y$, $x \mapsto c(x, y)$ is $\zeta$-semi-convex, there exists
          a $\max\left(\xi, \frac{\xi + \zeta}{\lambda}\right)$-strongly convex
          function $V : \X \to \Rsp$ such that the density of $\mu$ reads $f(x)
          = e^{-V(x)}$, and
        $$
        \alpha = 176\{ 1 + c_\infty \lambda^{-1} \}.
        $$
    \end{enumerate}
  \end{theorem}

  Under slightly more assumptions, the contraction constants $1 - \Theta(\lambda^2/c_\infty^2)$ shown under assumptions \ref{assump1} or \ref{assump2} can be shown to improve to $1 - \Theta(\lambda/c_\infty)$ after enough iterations have been performed. This is the content of the following corollary, proven in Section~\ref{sec:proof-corollary-large-t-contraction}.

  \begin{corollary}
  \label{cor:large-t-contraction}
      Assume that \ref{assump} holds and that either the cost $c$ is Lipschitz continuous and $\Y$ is bounded or the target measure $\nu$ is finitely supported. Then, there exists an integer $T$ depending on $\mu, \nu, c$ and $\psi_0$ such that for any $t \geq T$, guarantee \eqref{eq:exponential-convergence} of Theorem~\ref{thm:exponential-convergence} holds with a value of $\alpha$ improved to: 
      \begin{enumerate}
          \item  
          $ \alpha = 176 \{ 1 +
        (c_\infty + \frac{\xi}{2} R_\X^2)\kappa \lambda^{-1} \}
        $ under \ref{assump1}.
          \item 
          $
        \alpha = 176\{ 1 +
        c_{\infty}  \lambda^{-1}\}$ under \ref{assump2}.
      \end{enumerate}
  \end{corollary}

  Proofs of the Theorems~\ref{thm:quadratic-cost-log-concave} and \ref{thm:quadratic-cost-bounded-density} follow directly from the main theorem and its corollary stated above. 
  \begin{proof}[Proof of Theorem~\ref{thm:quadratic-cost-log-concave}]
    We apply Theorem~\ref{thm:exponential-convergence} under the conditions \ref{assump} and \ref{assump3}.
    The linear cost function $c(x,y) = -\sca{x}{y}$
    is $0$-semi-concave and $0$-semi-convex.
    Thus, $\xi=\zeta=0$.
    Hence, $\max\left(\xi, \frac{\xi + \zeta}{\lambda}\right) = 0$, and the condition on the density of $\mu$ reduces to log-concavity (i.e., $V$ is required to be convex).
    This concludes the proof of Theorem~\ref{thm:quadratic-cost-log-concave}.
    %\hfill\qed
  \end{proof}
   \begin{proof}[Proof of Theorem~\ref{thm:quadratic-cost-bounded-density}]
     We apply Theorem~\ref{thm:exponential-convergence} under the conditions \ref{assump} and \ref{assump1}. The linear cost $c(x,y) = -\sca{x}{y}$ is $0$-semi-concave, thus $\xi=0$.
     With $\xi=0$, Theorem~\ref{thm:exponential-convergence} gives
     $\alpha = 176(1 + \kappa c_{\infty}\lambda^{-1} + c_{\infty}^{2}\lambda^{-2})$,
     valid for any $t \geq 0$, proving the first part of Theorem~\ref{thm:quadratic-cost-bounded-density}.
     The second part follows from Corollary \ref{cor:large-t-contraction}, which ensures the existence of an integer $T$ such that for $t \geq T$, the constant $\alpha$ improves to $\alpha = 176(1 + \kappa c_{\infty}\lambda^{-1})$, which proves the second part of Theorem~\ref{thm:quadratic-cost-bounded-density}.
  \end{proof}

Before proving Theorem~\ref{thm:exponential-convergence} in the next sub-section, we make the following few remarks.

%\LC{add a small discussion (perhaps with bullet points) to unpack a bit. In particular I'd give examples or edge cases for (A2) and (A3) since they're a slightly technical. Also the asymptotic rate and the Gaussian rates could be added here. finally we could write 2-3 sentences on the proof technique: the key parts are Prop 4.4 and L2 strong convexity estimates for $\Kant$.}

\begin{remark}
    The convexity assumption made on the support $\X$ of $\mu$ is restrictive. 
    This assumption is needed in only one of the intermediate results, 
    Proposition~\ref{prop:strong-concavity} below, that proves a strong-concavity 
    estimate for the entropic Kantorovich functional. It was shown however 
    in \cite{Carlier2023} that such strong-concavity estimates could be obtained under 
    milder assumptions on $\X$, e.g. by assuming only that $\X$ is a 
    finite connected union of convex sets. Overall, we conjecture that assuming that 
    $\mu$ satisfies a Poincaré inequality should be a sufficient condition for 
    our results to hold.
    %\LC{finetune this paragraph a bit}
    %\tomas{My suggestion is to keep only the following:}
    %\tomas{The convexity assumption made on the support of $\mu$ can be relaxed; see \cite{Carlier2023} for details.}
\end{remark}

\begin{remark}
  The semi-concavity and semi-convexity assumptions made on $c$ are easily satisfied whenever 
  $c$ is smooth and $\X$ and $\Y$ are compact, since in this case $c$ is 
  both $\nr{c}_{\Class^2(\X \times \Y)}$-semi-concave and 
  $\nr{c}_{\Class^2(\X \times \Y)}$-semi-convex.
\end{remark}

\begin{remark}
\label{rk:approx-c-transform}
    Assumptions \ref{assump1}, \ref{assump2} or \ref{assump3} of Theorem~\ref{thm:exponential-convergence} all require the measure $\mu$ to be absolutely continuous. Because of this assumption, the $(c,\lambda)$-transform 
    $$ \phi^{c,\lambda}(y) = -\lambda \log \int e^{\frac{\phi(x) - c(x,y)}{\lambda}} \dd \mu(x)$$
    appearing in the definition of Sinkhorn's iterates \eqref{eq:sinkhorn-iterates} cannot be implemented exactly on a computer in general. In this sense, the convergence result of Theorem~\ref{thm:exponential-convergence} concerns an \emph{infinite dimensional} or \emph{continuous in space} version of Sinkhorn's algorithm. We notice however in Section~\ref{sec:convergence-rates-approximated-sequence} of the appendix that the results of Theorem~\ref{thm:exponential-convergence} can be extended to the practical setting where the above $(c, \lambda)$-transform is only approximated. Let us consider for instance that we access a discrete approximation $\hat{\mu} = \frac{1}{N} \sum_{i=1}^N \delta_{x_i}$ of $\mu$ that is such that for any potential $\phi \in L^\infty(\mu)$ having the same regularity as the cost function $c$,
    $$ \nr{\phi^{\widehat{c,\lambda}} - \phi^{c, \lambda}}_\infty \leq \eps$$
    for some fixed $\eps > 0$, where $\phi^{\widehat{c, \lambda}}(y) = -\lambda \log \frac{1}{N} \sum_{i=1}^N e^{\frac{\phi(x_i) - c(x_i,y)}{\lambda}}$. Then assuming \ref{assump} and either \ref{assump1} or \ref{assump2} for $\mu$ (no further assumption is made on $\hat{\mu}$), Proposition~\ref{prop:convergence-rates-approximated-sequence} ensures that that the \emph{approximated} Sinkhorn iterates $\hat{\psi}_{t+1} = (\hat{\psi}_t^{c,\lambda})^{\widehat{c,\lambda}}$ satisfy
    $$ E(\psi^{*}) - E(\hat{\psi}_{t}) \leq (1 - \alpha^{-1})^t(E(\psi^{*}) - E(\hat{\psi}_{0})) + 2\alpha \eps. $$
    This corresponds to an exponential contraction of the dual \emph{population} objective $E$ (that is associated to $\mu$), up to an error term which is proportional to the error made in the approximation of the $(c,\lambda)$-transforms.
\end{remark}

\begin{remark}
  The semi-concavity and semi-convexity constants $\xi, \zeta$ do not appear in the expressions 
  of $\alpha$ obtained in 
  Theorem~\ref{thm:exponential-convergence} under assumptions \ref{assump2}, \ref{assump3}. From the proofs, this fact may be interpreted with the idea that the strong-convexity of 
  the potential $V$ \emph{exactly compensates} the lack of concavity or convexity of the cost 
  $c$. Surprisingly, it does not seem that assuming further convexity for the potential $V$ (e.g. in \ref{assump2} assuming that $V$ is $\tilde{\xi}$-strongly convex for some $\tilde{\xi}>\xi$) can easily be leveraged to improve the contraction constants. 
\end{remark}

  \begin{remark}
      The statement of Corollary~\ref{cor:large-t-contraction} is asymptotic. 
      However, as detailed in Section~\ref{sec:refinement-cor-value-T} of the Appendix, the value 
      of the integer $T$ present in this statement can be made explicit in some cases. For instance, whenever the target measure $\nu$ is finitely supported over $\Y = \{y_1, \dots, y_N\}$, we have
      $$T = \Bigg\lceil\frac{\log\left(\frac{4 C_1 \delta_0}{\min_i \nu(y_i)^2}\right)}{-\log(1-\alpha^{-1})}\Bigg\rceil,$$
      where $\alpha = 176 \{ 1 +
        (c_\infty + \frac{\xi}{2} R_\X^2)\kappa \lambda^{-1}
        + c_\infty^2 \lambda^{-2} \}
        $ and $C_1 = 11\{
        (c_{\infty} + \frac{\xi}{2}R_{\mathcal{X}}^2) \kappa
        +
        \lambda
        +
        c_{\infty}^{2}\lambda^{-1}\}$ under \ref{assump1}; or $
        \alpha = 176\{ 1 +
        c_{\infty}  \lambda^{-1}
        +
        c_{\infty}^{2}\lambda^{-2}
        \}$ and 
        $C_1 = 11\{
        c_{\infty} 
        +
        \lambda
        +
        c_{\infty}^{2}\lambda^{-1}\}$ 
        under \ref{assump2}.
  \end{remark}

\subsection{Proof of Theorem~\ref{thm:exponential-convergence}}
\label{sec:exponential-convergence-thm-proof}

The proof of Theorem~\ref{thm:exponential-convergence} follows from the following two key independent propositions stated below.

The first proposition relates the sub-optimality gap $\delta_{t}$ to the improvement made in one step $\delta_{t} - \delta_{t+1}$. Prior works that prove polynomial bounds on $\delta_{t}$ (e.g., \cite{kalantari2008complexity, altschuler2017near, dvurechensky2018computational})
typically rely on bounds of the form
\begin{equation}
    \label{eq:hoeffding-type-bound}
    \delta_{t} \leq c_{\infty} \sqrt{2\lambda^{-1} (\delta_{t} - \delta_{t+1})},
\end{equation} from which the estimate $\delta_{t} = O(1/t)$ follows using traditional arguments in smooth convex optimization.
Instead of using the bound \eqref{eq:hoeffding-type-bound}, we prove a stronger result in the proposition below, which replaces the $c_{\infty}$ factor by a certain variance term related to the curvature of the semi-dual functional $E$. Our improvement comes from the use of Bernstein's inequality instead of Hoeffding's inequality to control cumulant--generating functions of bounded random variables appearing in the proof of Proposition~\ref{prop:one-step-improvement}; see Section~\ref{sec:proof-of-one-step-improvement-proposition} for details.
\begin{proposition}
  \label{prop:one-step-improvement}
  For any $t \geq 0$, the Sinkhorn iterates $(\psi_{s})_{s \geq 0}$ defined in
  \eqref{eq:sinkhorn-iterates} satisfy
  \begin{equation}
    \delta_{t}
    \leq
    2\sqrt{\lambda^{-1}
        \Var_{\nu}(\psi^{*} - \psi_{t})
      (\delta_{t} - \delta_{t+1})}
      + \frac{14c_{\infty}}{3}\lambda^{-1}(\delta_{t} - \delta_{t+1}).
  \end{equation}
\end{proposition}

\begin{remark}
    The proof of Proposition~\ref{prop:one-step-improvement} only requires the cost function $c$ to be of bounded oscillation norm. All the other assumptions enter through the strong-concavity estimates obtained in Proposition~\ref{prop:bound-variance-subopt} below.
\end{remark}

The second proposition relates $\Var_{\nu}(\psi^{*} - \psi_{t})$ back to $\delta_{t}$ and $\delta_{t+1}$ following an estimation of the strong-concavity of the semi-dual functional $E$. The proof of this proposition, given in Section~\ref{sec:proof-of-variance-subopt-bound}, essentially relies on the Prékopa-Leindler inequality. It extends to bounded semi-concave costs the existing estimations of the strong-concavity of the entropic and non-regularized Kantorovich functional that were previously proven for a linear cost in \cite{delalande2022concavity-sdeot, delalande2021concavity-quad-ot}.

\begin{proposition}
  \label{prop:bound-variance-subopt} Assume \ref{assump}. Then for any $t \geq 0$, the Sinkhorn iterates $(\psi_{s})_{s \geq 0}$ defined in \eqref{eq:sinkhorn-iterates} satisfy
    $$
      \Var_{\nu}(\psi^{*} - \psi_{t})
      \leq
      C_1 \delta_t + C_2(\delta_t - \delta_{t+1})
    $$
in either of the following cases:
\begin{enumerate}
    \item \ref{assump1} holds, $C_1 = 11\{
        (c_{\infty} + \frac{\xi}{2}R_{\mathcal{X}}^2) \kappa
        +
        \lambda
        +
        c_{\infty}^{2}\lambda^{-1}\}$ and $C_2 = 0$.
      \item \ref{assump2} holds, $C_1 = 11\{
        c_{\infty}
        +
        \lambda
        +
        c_{\infty}^{2}\lambda^{-1}\}$ and $C_2 = 0$.
      \item \ref{assump3} holds, $C_1 = 11 \left( c_{\infty} + \lambda \right)$ and $C_2 = 3 c_{\infty}^2 \lambda^{-1}$.
     % \item \ref{assump1} hold, $t \geq T$ for large enough
     %   $T=T(\mu,\nu,c,\psi_{0})$,
     % $C_{1} = 6\lambda + 6e(c_{\infty} +
     % \frac{\xi}{2}R_{\mathcal{X}}^{2})\kappa$ and $C_2 = 0$. \tomas{added a proof sketch for this.}
\end{enumerate}
\end{proposition}
Combining Proposition~\ref{prop:one-step-improvement} with
Proposition~\ref{prop:bound-variance-subopt} yields convergence guarantees of
the form $\delta_{t+1} \leq C'\delta_{t}$ for some contraction rate $C' \in
(0,1)$. This is shown in the following lemma.
\begin{lemma}
    \label{lemma:deriving-contraction-rate}
    Suppose that $\|c\|_{\mathrm{osc}} = c_{\infty} < \infty$. Let $(\psi_{s})_{s \geq 0}$ be the Sinkhorn iterates.
    If for some $t$ we have
    $
      \Var_{\nu}(\psi^{*} - \psi_{t})
      \leq
      C_1 \delta_t + C_2(\delta_t - \delta_{t+1})
    $
    then
    $$
      \delta_{t+1} \leq (1-\alpha^{-1})\delta_{t},\quad\text{where}\quad
      \alpha =
      \max\left\{
        16\lambda^{-1}C_{1},
        4\sqrt{\lambda^{-1}C_{2}}
          +
          \frac{28c_{\infty}}{3}\lambda^{-1}
        \right\}.
    $$
\end{lemma}
The proof of the above lemma involves straightforward algebraic manipulations;
we defer the details to
Appendix~\ref{sec:proof-of-deriving-contraction-rate-lemma}. We will now
complete the proof of Theorem~\ref{thm:exponential-convergence}.
\begin{enumerate}
  \item If \ref{assump1} holds, then
    Propositions~\ref{prop:one-step-improvement} and \ref{prop:bound-variance-subopt} with
    Lemma~\ref{lemma:deriving-contraction-rate} yield $\delta_{t+1} \leq (1-\alpha^{-1})\delta_{t}$ with
    \begin{align}
        \alpha &= \max\left( 176 \lambda^{-1}\left\{ (c_\infty + \frac{\xi}{2} R_\X^2)\kappa + \lambda + c_\infty^2 \lambda^{-1} \right\}, \frac{28 c_\infty}{3} \lambda^{-1} \right)  \\
        &= 176 \left\{ 1 + (c_\infty + \frac{\xi}{2} R_\X^2)\kappa \lambda^{-1} + c_\infty^2 \lambda^{-2} \right\},
    \end{align}
    where we used the fact that $\kappa \geq 1$.
  \item If \ref{assump2} holds, then
    Propositions~\ref{prop:one-step-improvement} and \ref{prop:bound-variance-subopt} with
    Lemma~\ref{lemma:deriving-contraction-rate} yield $\delta_{t+1} \leq (1-\alpha^{-1})\delta_{t}$ with
    \begin{align}
        \alpha &= \max\left( 176 \lambda^{-1}\left\{ c_\infty + \lambda + c_\infty^2 \lambda^{-1} \right\}, \frac{28 c_\infty}{3} \lambda^{-1} \right)  \\
        &= 176 \left\{ 1 + c_\infty \lambda^{-1} + c_\infty^2 \lambda^{-2} \right\}.
    \end{align}
  \item If \ref{assump3} holds, then
    Propositions~\ref{prop:one-step-improvement} and
    \ref{prop:bound-variance-subopt} with
    Lemma~\ref{lemma:deriving-contraction-rate} yield $\delta_{t+1} \leq (1-\alpha^{-1})\delta_{t}$ with
    \begin{align}
        \alpha &= \max\left( 176 \lambda^{-1}\left\{ c_\infty + \lambda \right\}, 4\sqrt{3 c_\infty^2 \lambda^{-2}} + \frac{28 c_\infty}{3} \lambda^{-1} \right)  \\
        &= 176 \left\{ 1 + c_\infty \lambda^{-1} \right\}.
    \end{align}
  %\item If \ref{assump1} holds and $t \geq T$ for large enough $T$, then
  % Propositions~\ref{prop:one-step-improvement} and
  % \ref{prop:bound-variance-subopt} with
  % Lemma~\ref{lemma:deriving-contraction-rate} yield $\delta_{t+1} \leq (1-\alpha^{-1})\delta_{t}$ with
  % \begin{align}
  %     \alpha
  %     &=
  %     \max\left(
  %       96
  %       +
  %       96e\kappa\frac{c_{\infty} + \frac{\xi}{2}R_{\mathcal{X}}^{2}}{\lambda}
  %       ,
  %       \frac{28 c_\infty}{3} \lambda^{-1}
  %      \right)
  %      \\
  %     &=
  %     96
  %     +
  %     96e\kappa (c_{\infty} + \frac{\xi}{2}R_{\mathcal{X}}^{2}) \lambda^{-1},
  % \end{align}
  % where we used the fact that $\kappa \geq 1$.
\end{enumerate}  \hfill\qed

\section{One-step improvement bound (Proof of Proposition~\ref{prop:one-step-improvement})}
\label{sec:proof-of-one-step-improvement-proposition}

By the concavity of $\Kant$ guaranteed in Lemma~\ref{lemma:K-concave-derivatives},
the semi-dual functional is concave (this also follows directly from the definition of $E$ as a partial supremum of concave functions). Hence, the first derivative computed in
Lemma~\ref{lemma:K-concave-derivatives} with $\psi = \psi_{t}$ and $v =
\psi^{*} - \psi_{t}$ yields
\begin{equation}
  \delta_{t} \leq
  \sca{\nu - \nu[\psi_{t}]}{\psi^{*} - \psi_{t}}.
\end{equation}
In particular, for any $\eta > 0$ we have
\begin{align}
  \delta_{t}
  &\leq
  \eta^{-1}
  \left\{
    \sca{\nu - \nu[\psi_{t}]}{\eta(\psi^{*} - \psi_{t})}
    - \kl{\nu}{\nu[\psi_{t}]}
  \right\}
  +
  \eta^{-1}
  \kl{\nu}{\nu[\psi_{t}]}
  \\
  &\leq
  \eta^{-1}\sup_{\nu' \in \mathcal{P}(\mathcal{X})}
  \left\{
    \sca{\nu' - \nu[\psi_{t}]}{\eta(\psi^{*} - \psi_{t})}
    - \kl{\nu'}{\nu[\psi_{t}]}
  \right\}
  +
  \eta^{-1}
  \kl{\nu}{\nu[\psi_{t}]}
  \\
  &=\eta^{-1}  \log \Esp_{ \nu[\psi_{t}]}\left[
    \exp\left(
      \eta f
    \right)
  \right]
  +
  \eta^{-1}
  \kl{\nu}{\nu[\psi_{t}]},
  \label{eq:delta-t-bound-middle-step}
\end{align}
where $f(y)= \psi^{*} - \psi_{t} - \Esp_{\nu[\psi_{t}]}[\psi^{*} - \psi_{t}]$. The last line above follows from the expression of the convex conjugate of $\kl{\cdot}{\nu[\psi_{t}]}$ and 
$$\sca{\nu' - \nu[\psi_{t}]}{\eta(\psi^{*} - \psi_{t})}=\sca{\nu'}{\eta(\psi^{*} - \psi_{t})-\eta\sca{\nu[\psi_{t}]}{\psi^{*} - \psi_{t})}}=\sca{\nu'}{\eta f}.$$

By Bernstein's moment generating function bound for bounded random variables~\cite[Proposition~2.14]{wainwright2019high} and the fact
that $\psi^{*} - \psi_{t}$ and hence $f$ is contained in an interval of length at most $2c_{\infty}$, it follows that for any $\eta \in (0, \frac{3}{2c_{\infty}})$
we have
\begin{align}
  \log \Esp_{y\sim \nu[\psi_{t}]}\left[
    \exp\left(
      \eta f
    \right)
  \right]
  \leq
  \frac{
    \eta^{2}\Var_{\nu[\psi_{t}]}(\psi^{*} - \psi_{t})
  }{
    2(1 - \eta\frac{2c_{\infty}}{3})
  }.
\end{align}
Since this holds for any $\eta \in (0,
\frac{3}{2c_{\infty}})$, we obtain from bound \eqref{eq:delta-t-bound-middle-step} the estimate
\begin{align}
  \delta_{t}
  &\leq
  \inf_{0 < \eta < \frac{3}{2c_{\infty}}}
  \left\{
    \frac{
      \eta
      \Var_{\nu[\psi_{t}]}(\psi^{*} - \psi_{t})
    }{
      2(1 - \eta\frac{2c_{\infty}}{3})
    }
    +
    \eta^{-1}\kl{\nu}{\nu[\psi_{t}]}
  \right\}
  \\
  &\leq
  \sqrt{2
      \Var_{\nu[\psi_{t}]}(\psi^{*} - \psi_{t})
      \kl{\nu}{\nu[\psi_{t}]}}
  + \frac{2c_{\infty}}{3}\kl{\nu}{\nu[\psi_{t}]},
  \label{eq:final-corollary-from-bernsteins-bound}
\end{align}
where the final inequality follows by optimizing in $\eta$ (see details in~\cite{boucheron2013concentration}, Lemma~2.4 and the computation following Equation~(2.5) on page 29 therein).

By equation~\eqref{eq:sinkhorn-iterates-improvement}, we have
\begin{equation}
    \kl{\nu}{\nu[\psi_t]} \leq
     \lambda^{-1}(\delta_t - \delta_{t+1})
     \label{eq:kl-one-step-improvement-bound}
\end{equation}
Combining the above inequality with
\eqref{eq:final-corollary-from-bernsteins-bound} yields
\begin{align}
  \delta_{t}
  \leq
  \sqrt{2\lambda^{-1}
      \Var_{\nu[\psi_{t}]}(\psi^{*} - \psi_{t})
    (\delta_{t} - \delta_{t+1})}
    + \frac{2c_{\infty}}{3}\lambda^{-1}(\delta_{t} - \delta_{t+1}).
    \label{eq:one-step-improvement-before-replacing-variance}
\end{align}
In order to complete the proof of Proposition~\ref{prop:one-step-improvement}, it remains to replace the $\Var_{\nu[\psi_{t}]}(\psi^{*} - \psi_{t})$ term in the above inequality with $\Var_{\nu}(\psi^{*} - \psi_{t})$.
To this end, we will use the following lemma proved in Appendix~\ref{sec:proof-of-variance-comparison-inequality}.
\begin{lemma}
  \label{lemma:variance-comparison-inequality}
  Let $\rho,\pi \in \mathcal{P}(\mathcal{X})$
  by any probability measures such that $\rho \ll \pi$ and $\pi \ll \rho$.
  Then, for any function $f : \mathcal{X} \to [a,b]$ it holds that
  \begin{equation}
    \Var_{\rho}(f)
    \geq
    \frac{1}{2}\Var_{\pi}(f)
    - (b-a)^2 \mathrm{min}(
      \kl{\rho}{\pi}, \kl{\pi}{\rho}
    ).
  \end{equation}
\end{lemma}

Indeed, by Lemma~\ref{lemma:variance-comparison-inequality} and the fact that $\|\psi^{*} - \psi_{t}\|_{\mathrm{osc}} \leq 2c_{\infty}$, we have
\begin{align}
    \Var_{\nu[\psi_{t}]}(\psi^{*} - \psi_{t})
    &\leq 2\Var_{\nu}(\psi^{*} - \psi_{t})
    + 8 c_{\infty}^2 \kl{\nu}{\nu[\psi_{t}]}
    \\
    &\leq
    2\Var_{\nu}(\psi^{*} - \psi_{t})
    + 8 c_{\infty}^2 \lambda^{-1}(\delta_{t} - \delta_{t+1}),
\end{align}
where the second inequality follows from \eqref{eq:kl-one-step-improvement-bound}.
Plugging the above inequality into \eqref{eq:one-step-improvement-before-replacing-variance}
and applying the bound $\sqrt{a + b} \leq \sqrt{a} + \sqrt{b}$, valid for $a,b > 0$, we obtain
\begin{align}
  \delta_{t}
  \leq
  2\sqrt{\lambda^{-1}
      \Var_{\nu}(\psi^{*} - \psi_{t})
    (\delta_{t} - \delta_{t+1})}
    + \frac{14c_{\infty}}{3}\lambda^{-1}(\delta_{t} - \delta_{t+1}),
\end{align}
which completes the proof of Proposition \ref{prop:one-step-improvement}. \hfill\qed

\section{Strong-concavity estimate (Proof of Proposition~\ref{prop:bound-variance-subopt})}
\label{sec:proof-of-variance-subopt-bound}

\subsection{Proof of Proposition~\ref{prop:bound-variance-subopt}}

By optimality of $\psi^*$, one has $\nu[\psi^*] = \nu.$ Thus, for $v = \psi^* - \psi_t$, one has from Lemma~\ref{lemma:K-concave-derivatives} the following expression for the dual sub-optimality gap:
\begin{align} \label{eq:reasonning}
    \delta_t &= \Kant(\psi^*) - \Kant(\psi_t) - \frac{\dd}{\dd \eps} \Kant(\psi_t + \eps v)\bigg|_{\eps = 1} \notag \\
    &= \int_{\eps=0}^1 \frac{\dd}{\dd \eps} \Kant(\psi_t + \eps v) \dd \eps - \frac{\dd}{\dd \eps} \Kant(\psi_t + \eps v)\bigg|_{\eps = 1} \notag  \\
    &=  - \int_{\eps=0}^1 \int_{s=\eps}^{1} \frac{\dd^2}{\dd s^2} \Kant(\psi_t + s v) \dd s \dd \eps.
\end{align}
With these computations, we see that proving Proposition~\ref{prop:bound-variance-subopt} reduces to proving a \emph{strong-concavity estimate} for $\Kant$ of the type
\begin{equation} \label{eq:wanted-strong-concavity-estimate}
  \frac{\dd^2}{\dd s^2} \Kant(\psi_t + s v) \leq -C \Var_{\nu[\psi_t + sv]}(v),
\end{equation}
for some constant $C$. Associativity of variances ensures
     \begin{equation*}
         \Var_{\nu[\psi_t + sv]}(v) = \Var_{x \sim \mu}(\Esp_{\nu_x[\psi_t + sv]}(v)) + \Esp_{x \sim \mu} \Var_{\nu_x[\psi_t + sv]}(v),
     \end{equation*}
     and Lemma~\ref{lemma:K-concave-derivatives} guarantees that $\frac{\dd^2}{\dd s^2} \Kant(\psi_t + s v) = - \frac{1}{\lambda} \Esp_{x \sim \mu} \Var_{\nu_x[\psi_t + sv]}(v)$. We thus get directly
     \begin{equation} \label{eq:associativity-var}
         \frac{\dd^2}{\dd s^2} \Kant(\psi_t + s v) = -\frac{1}{\lambda} ( \Var_{\nu[\psi_t + sv]}(v) - \Var_{x \sim \mu}(\Esp_{\nu_x[\psi_t + sv]}(v)) ).
     \end{equation}
In order to conclude and get a strong-concavity estimate of the form of \eqref{eq:wanted-strong-concavity-estimate}, there remains to compare $\Var_{x \sim \mu}(\Esp_{\nu_x[\psi_t + sv]}(v))$ to either $\frac{\dd^2}{\dd s^2} \Kant(\psi_t + s v) $ or $\Var_{\nu[\psi_t + sv]}(v) $:
this is the content of the following proposition, which is proven in Section~\ref{sec:strong-concavity} and might be of independent interest.
\begin{proposition} \label{prop:strong-concavity} Assume that \ref{assump} holds. Then, for any $\psi, v : \Rsp^d \to \Rsp$ and $s \in \Rsp$,
the \emph{entropic Kantorovich functional} satisfies
    \begin{equation*}
        \frac{\dd^2}{\dd s^2} \Kant(\psi + s v)
        \leq
        -\left(
          C
            +  \lambda
        \right)^{-1}
        \Var_{\nu[\psi + s v]}(v)
    \end{equation*}
    in either of the following cases:
    \begin{enumerate}
        \item \ref{assump1} holds and $C = e(c_{\infty} + \frac{\xi}{2}R_{\mathcal{X}}^2) \kappa$.
        \item \ref{assump2} or \ref{assump3} holds and $C = e c_\infty$.
    \end{enumerate}
\end{proposition}
\begin{remark}
    The above strong-concavity estimate and its proof idea relying on the Prékopa-Leindler inequality already appeared in \cite[Theorem 3.2]{delalande2022concavity-sdeot}, but only in the case where $c$ was the linear ground cost $c(x,y) = -\sca{x}{y}$. Interestingly, letting the regularization strength $\lambda$ go to $0$ in the above estimate allows to also get an estimation of the strong-concavity of the \emph{non-regularized} Kantorovich functional, which is the functional
    $$ \Kant : \psi \mapsto \int \psi^c \dd \mu,$$
    where $\psi^c(\cdot) = \inf_y c(\cdot,y) - \psi(y)$ is the usual $c$-transform. This allows to extend to bounded semi-concave costs the strong-concavity estimate that was already established in \cite[Theorem 2.1]{delalande2021concavity-quad-ot} for the non-regularized Kantorovich functional with linear ground cost $c(x,y) = -\sca{x}{y}$, which in turn was proven relying on the Brascamp-Lieb concentration inequality.
\end{remark}

Before proving Proposition~\ref{prop:strong-concavity} in the next sub-section, we conclude the proof of Proposition~\ref{prop:bound-variance-subopt}. From the reasoning in \eqref{eq:reasonning}, we know that the sub-optimality gap admits the following expression:
    \begin{equation*}
        \delta_t =  - \int_{\eps=0}^1 \int_{s=\eps}^{1} \frac{\dd^2}{\dd s^2} \Kant(\psi_t + s v) \dd s \dd \eps.
    \end{equation*}
    Proposition~\ref{prop:strong-concavity} then ensures the following lower bound:
    \begin{align}
      \label{eq:lower-bound-delta_t-var}
      \left( C + \lambda \right)
      \delta_t
      \geq
      \int_{\eps=0}^1 \int_{s=\eps}^{1} \Var_{\nu[\psi_t + s v]}(v) \dd s \dd \eps,
    \end{align}
    where $C = e(c_{\infty} + \frac{\xi}{2}R_{\mathcal{X}}^2) \kappa$ under \ref{assump1} and $C = ec_\infty$ under \ref{assump2} or \ref{assump3}.
    Recall that $v = \psi^{*} - \psi_{t}$. Thus, $\nu[\psi_{t} + v] =
    \nu[\psi^{*}] = \nu$, and the oscillation norm of $v$ is bounded by
    $2 c_{\infty}$. It follows by
    Lemma~\ref{lemma:variance-comparison-inequality} that
    \begin{align}
      \left( C + \lambda
        \right)
      \delta_t
      &\geq
      \int_{\eps=0}^1 \int_{s=\eps}^{1}
      \frac{1}{2}\Var_{\nu}(v)
      - (2c_{\infty})^2 \kl{\nu}{\nu[\psi_t + sv]}
      \dd s \dd \eps
      \\
      &=
      \frac{1}{4}\Var_{\nu}(v)
      -
      4c_{\infty}^2
      \int_{\eps=0}^1 \int_{s=\eps}^{1}
      \kl{\nu}{\nu[\psi_t + sv]}\dd s \dd \eps.
      \label{eq:variance-delta-comparison-with-kl-terms}
    \end{align}
    Next, %recall from \eqref{eq:sinkhorn-iterates-improvement} that
    %$
    %  \kl{\nu}{\nu[\psi_{t}]} \leq \lambda^{-1}(\delta_{t} - \delta_{t+1})
    %  \leq \lambda^{-1}\delta_{t}.
    %$
    %Similarly,
    notice that for any $s \in [0,1]$ we have by equation~\eqref{eq:KL-nu-nu-psi}:
    \begin{align}
      \kl{\nu}{\nu[\psi_{t} + sv]} &= \frac{1}{\lambda}  \sca{ (\psi_t + sv)^{\overline{c,\lambda}} - (\psi_t + s v)}{\nu} \\
      &= \frac{1}{\lambda} \left(F((\psi_t + sv)^{c,\lambda}, (\psi_t + sv)^{\overline{c,\lambda}}) - F((\psi_t + sv)^{c,\lambda}, \psi_t + sv)\right) \\
      &\leq \frac{1}{\lambda} \left(F((\psi^*)^{c,\lambda}, \psi^*) - F((\psi_t + sv)^{c,\lambda}, \psi_t + sv)\right) \\
      &= \frac{1}{\lambda}\left(E(\psi^*) - E((1-s)\psi_t + s \psi^*)\right) \\
      &\leq \frac{1-s}{\lambda}\left(E(\psi^*) - E(\psi_t) \right) \\
      &= \frac{1-s}{\lambda} \delta_t,
      \label{eq:kl-s-bound}
    \end{align}
where the first inequality follows from the definition of $\psi^*$ as a maximizer for the semi-dual functional and the last inequality follows from the concavity of $E$.
\begin{remark}
The above inequality is where we lose a factor $\lambda/c_\infty$ in the final convergence rate under \ref{assump1}. Whenever this bound can be improved to $\lambda^{-1}(\delta_{t} -
      \delta_{t+1})$, the resulting convergence rate becomes  $(1 -\Theta(\lambda/c_{\infty}))^t$.
\end{remark}
Plugging the inequality \eqref{eq:kl-s-bound} in \eqref{eq:variance-delta-comparison-with-kl-terms} yields
\begin{align}
    \label{eq:upper-bound-var-subopt}
      \Var_{\nu}(v) &\leq
      \left(
        4C
        +
        4 \lambda
        +
        \frac{8}{3} c_{\infty}^2 \lambda^{-1}
      \right)\delta_{t}.
 \end{align}
    %\LC{I think that the $6$ in the first line should be $(32\times 4)/6\leq 22$. If so, many numerical constants need to be adjusted in the paper.}
    %\tomas{I changed the above constants up to here. The below will be wrong.}
    %\textcolor{green}{Alex (personal note): revoir constantes à partir d'ici}

    Under assumption \ref{assump1}, $C = e(c_{\infty} + \frac{\xi}{2}R_{\mathcal{X}}^2) \kappa$ so that the last bound entails
    \begin{align}
      \Var_{\nu}(v) \leq
      11\left(
       (c_{\infty} + \frac{\xi}{2}R_{\mathcal{X}}^2)\kappa
        +
        \lambda
        +
        c_{\infty}^{2}\lambda^{-1}
      \right)
      \delta_{t},
    \end{align}
    which completes the proof of the first bound claimed in Proposition~\ref{prop:bound-variance-subopt}.

    Under assumption \ref{assump2}, $C = ec_\infty$ so that bound \eqref{eq:upper-bound-var-subopt}
    becomes
    \begin{align}
      \Var_{\nu}(v) \leq
      11\left(
       c_{\infty}
        +
        \lambda
        +
        c_{\infty}^{2}\lambda^{-1}
      \right)
      \delta_{t},
    \end{align}
    which completes the proof of the second bound claimed in Proposition~\ref{prop:bound-variance-subopt}.

    Assuming now \ref{assump3}, one has an improved bound for $\kl{\nu}{\nu[\psi_{t} + sv]}$. Indeed, the following proposition (proved in Appendix~\ref{sec:convexity-sinkhorn-transformation}) ensures the pointwise convexity of the $(\overline{c, \lambda})$-transformation under assumption \ref{assump3}:
    \begin{proposition} \label{prop:convexity-sinkhorn-transformation}
        Under assumptions \ref{assump} and \ref{assump3}, for any $\psi_0, \psi_1 : \Rsp^d \to \Rsp$ and $\alpha \in [0,1]$, it holds
        $$ \left( (1-\alpha) \psi_0 + \alpha \psi_1 \right)^{\overline{c, \lambda}} \leq (1-\alpha) \left(  \psi_0 \right)^{\overline{c, \lambda}} + \alpha \left(  \psi_1 \right)^{\overline{c, \lambda}}.$$
    \end{proposition}
    Thus, for any $y \in \Rsp^d$, one has
    \begin{align*}
        (\psi_t + sv)^{\overline{c, \lambda}}(y) &= ((1-s)\psi_t + s\psi^*)^{\overline{c, \lambda}}(y) \\
        &\leq (1-s) \psi_t^{\overline{c, \lambda}}(y) + s (\psi^*)^{\overline{c, \lambda}}(y) \\
        &= (1-s) \psi_t^{\overline{c, \lambda}}(y) + s \psi^*(y),
    \end{align*}
    where $(\psi^*)^{\overline{c, \lambda}} = \psi^*$ follows from the optimality of $\psi^*$. This ensures the estimate
    \begin{align*}
        \kl{\nu}{\nu[\psi_{t} + sv]} &= \frac{1}{\lambda}  \sca{ (\psi_t + sv)^{\overline{c,\lambda}} - (\psi_t + s v)}{\nu} \\
        &\leq \frac{1-s}{\lambda} \sca{\psi_t^{\overline{c,\lambda}} - \psi_t}{\nu} \\
        &= (1-s)  \KL(\nu | \nu[\psi_t]) \\
        &\leq (1-s) \lambda^{-1}(\delta_{t} - \delta_{t+1}),
    \end{align*}
    where we used \eqref{eq:sinkhorn-iterates-improvement}. Plugging this bound in \eqref{eq:variance-delta-comparison-with-kl-terms} and using $C= ec_\infty$ yields
    \begin{align}
      \Var_{\nu}(v)
      &\leq
      4\left( ec_\infty + \lambda \right)
      \delta_t
      + \frac{8}{3} c_{\infty}^2 \lambda^{-1}(\delta_{t} - \delta_{t+1}) \\
      &\leq 11 \left( c_{\infty} + \lambda \right)
      \delta_t
      +  3 c_{\infty}^2 \lambda^{-1}(\delta_{t} - \delta_{t+1}).
    \end{align}
    This complete the proof of the last bound claimed in Proposition \ref{prop:bound-variance-subopt}.   \hfill\qed

%\tomas{For the non-asymptotic part}

%\LC{Perhaps same issue as above with the factor $38$.}

\subsection{Proof of Proposition~\ref{prop:strong-concavity}}
\label{sec:strong-concavity}

In order to prove this result, we give ourselves a positive Radon measure $\rho \in \mathcal{M}_+(\Rsp^d)$ and
introduce an associated functional $I_\rho : (\Rsp^d \to \Rsp) \to \Rsp$, defined by
$$ I_\rho : \psi \mapsto \log \int_\X \exp \left( \psi^{c, \lambda}(x) \right) \dd \rho(x).$$
The idea of this definition comes from~\cite[Exercise 2.2.11]{klartag-lecture-notes}. From Lemma~\ref{lemma:derivatives-c-transform}, the mapping $s \mapsto I_\rho(\psi + sv)$ is $\Class^2$ on $\Rsp$. Denote $\tilde{\mu}[\psi]$ the probability measure in $\Prob(\X)$ satisfying
$$ \dd \tilde{\mu}[\psi](x) = \frac{\exp \left( \psi^{c, \lambda}(x) \right)\big|_\X}{\int_\X \exp \left( \psi^{c, \lambda}(\tilde{x}) \right) \dd \rho(\tilde{x})} \dd \rho(x).$$
Then the derivatives of $I$ are the following:
\begin{align}
    \frac{\dd}{\dd s} I_\rho(\psi + s v) &= \sca{ \frac{\dd}{\dd s} (\psi + s v)^{c ,\lambda} }{ \tilde{\mu}[\psi + s v] }, \label{eq:der-I}\\
    \frac{\dd^2}{\dd s^2} I_\rho(\psi + s v) &= \sca{ \frac{\dd^2}{\dd s^2} (\psi + s v)^{c ,\lambda} }{ \tilde{\mu}[\psi + s v] } + \Var_{ \tilde{\mu}[\psi + s v] }( \frac{\dd}{\dd s} (\psi + s v)^{c ,\lambda} ). \label{eq:der-2-I}
\end{align}
Proposition~\ref{prop:strong-concavity} will follow from the above expression of the second order derivative of $I_\rho$ and from the next result that ensures the concavity of $I_\rho$ for specific $\rho$ as a consequence of the Prékopa-Leindler inequality \cite{prekopa1, leindler, prekopa2}. Note that the following lemma is the only result that uses (and necessitates) the convexity assumption made on $\X$ in \ref{assump}.
\begin{lemma} \label{lemma:I-concave}
    If assumption \ref{assump} holds and $\rho = e^{-W}$ where $W : \Rsp^d \to \Rsp$ is a $\xi$-strongly convex function, then the functional $I_\rho$ is concave.
\end{lemma}
The above lemma is proved at the end of this section. We will now complete the
proof of Proposition~\ref{prop:strong-concavity}.

Under assumption \ref{assump1}, we choose $\rho = e^{-\frac{\xi}{2}\nr{x}^2}$. Otherwise, under assumptions
\ref{assump2} or \ref{assump3}, we choose $\rho = \mu$. In any of these cases, $\rho = e^{-W}$ with $W : \Rsp^d \to \Rsp$ a  $\xi$-strongly convex function.
We are thus in position to apply Lemma~\ref{lemma:I-concave}, that ensures that $I_\rho$ is a concave functional. From this concavity, we deduce
    $$ \frac{\dd^2}{\dd s^2} I_\rho(\psi + s v) \leq 0,$$
    that is with formula \eqref{eq:der-2-I}:
    \begin{equation} \label{eq:strong-concavity-1}
        \sca{ \frac{\dd^2}{\dd s^2} (\psi + s v)^{c ,\lambda} }{ \tilde{\mu}[\psi + s v] } \leq - \Var_{ \tilde{\mu}[\psi + s v] }( \frac{\dd}{\dd s} (\psi + s v)^{c ,\lambda} ) \leq 0.
    \end{equation}
    In order to conclude, we will compare the density of $\tilde{\mu}[\psi + s v]$ to the density of $\mu$. We have that
    $$ \dd \tilde{\mu}[\psi+ s v](x) = \frac{\exp \left( (\psi + sv)^{c, \lambda}(x) \right)\big|_\X}{\int_\X \exp \left( (\psi + sv)^{c, \lambda}(\tilde{x}) \right) \dd \rho(\tilde{x})} \dd \rho(x) = \frac{\exp \left( (\psi + sv)^{c, \lambda}(x) \right)\big|_\X}{Z} \dd \rho(x).$$
    From Lemma~\ref{lemma:properties-c-transforms}, we know that $(c,\lambda)$-transforms have finite oscillation upper bounded by $c_{\infty}$. In particular, there exists finite constants $m, M \in \Rsp$ such that $m \leq M$, $M - m \leq c_\infty$ and
    $$ m \leq (\psi + sv)^{c, \lambda} \leq M. $$
    This gives the bounds
    \begin{equation}
        \label{eq:comparison-tilde-mu-rho}
        \frac{e^{m}}{Z} \dd \rho(x) \leq \dd \tilde{\mu}[\psi+ s v](x) \leq \frac{e^M}{Z}\dd \rho(x).
    \end{equation}
    We finalize the comparison of $\tilde{\mu}$ to $\mu$ using the specific choices of $\rho$ we made under the different assumptions.

    \noindent \textbf{Bound under assumption \ref{assump1}.} Under assumption \ref{assump1}, we chose $\rho = e^{-\frac{\xi}{2}\nr{x}^2}$. But under this assumption, $\X \subset B(0, R_\X)$, we thus have for $x \in \X$ the bound
    $$ 0 \leq \nr{x} \leq R_\X.$$
    Injecting this bound into \eqref{eq:comparison-tilde-mu-rho} gives under \ref{assump1} the following bound on $\X$:
    $$ \frac{e^{m-\frac{\xi}{2} R_\X^2}}{Z}  \leq \tilde{\mu}[\psi+ s v] \leq \frac{e^M}{Z}.$$
    Now recall that under \ref{assump1}, the density $f$ of $\mu$ over $\X$ is upper- and lower-bounded. We denote under this assumption $m_\mu = \inf_{x \in \X} f(x) > 0$ and  $M_\mu = \sup_{x \in \X} f(x) <+\infty$, so that $\kappa = \frac{M_\mu}{m_\mu}$. Combining the previous estimate with the bound
    $$ m_\mu \leq \mu \leq M_\mu $$
    satisfied by $\mu$ on $\X$, we get the comparison
    $$ \frac{e^{m-\frac{\xi}{2} R_\X^2}}{Z M_\mu} \mu \leq  \tilde{\mu}[\psi+ s v] \leq \frac{e^{M}}{Z m_\mu} \mu. $$
    Inequality \eqref{eq:strong-concavity-1} thus ensures under assumption \ref{assump1}
    \begin{equation*}
       \frac{e^M}{Z m_\mu} \sca{ \frac{\dd^2}{\dd s^2} (\psi + s v)^{c ,\lambda} }{ \mu } \leq - \frac{e^{m-\frac{\xi}{2} R_\X^2}}{Z M_\mu} \Var_{ \mu }( \frac{\dd}{\dd s} (\psi + s v)^{c ,\lambda} ),
    \end{equation*}
which gives, under assumption \ref{assump1},
\begin{equation}
    \label{eq:strong-concavity-4-assump1}
    \kappa e^{c_\infty +\frac{\xi}{2}R_\X^2}  \sca{ \frac{\dd^2}{\dd s^2} (\psi + s v)^{c ,\lambda} }{ \mu } \leq - \Var_{ \mu }( \frac{\dd}{\dd s} (\psi + s v)^{c ,\lambda} ).
\end{equation}

\noindent \textbf{Bound under assumption \ref{assump2} or \ref{assump3}.} Now, if we no longer assume \ref{assump1} but instead assume \ref{assump2} or \ref{assump3}, the choice $\rho = \mu$ made above directly implies in \eqref{eq:comparison-tilde-mu-rho} the comparison
$$ \frac{e^{m}}{Z } \mu \leq  \tilde{\mu}[\psi+ s v] \leq \frac{e^{M}}{Z} \mu. $$
Inequality \eqref{eq:strong-concavity-1} thus ensures under assumptions \ref{assump2} or \ref{assump3}
\begin{equation}
    \label{eq:strong-concavity-4-assump-2-3}
    e^{c_\infty}  \sca{ \frac{\dd^2}{\dd s^2} (\psi + s v)^{c ,\lambda} }{ \mu } \leq - \Var_{ \mu }( \frac{\dd}{\dd s} (\psi + s v)^{c ,\lambda} ).
\end{equation}

\noindent \textbf{Concluding the proof of Proposition~\ref{prop:strong-concavity}.} By inequalities \eqref{eq:strong-concavity-4-assump1} and \eqref{eq:strong-concavity-4-assump-2-3},
we thus have in any case a bound of the type
\begin{equation}
    \label{eq:strong-concavity-4}
    C  \sca{ \frac{\dd^2}{\dd s^2} (\psi + s v)^{c ,\lambda} }{ \mu } \leq - \Var_{ \mu }( \frac{\dd}{\dd s} (\psi + s v)^{c ,\lambda} ),
\end{equation}
where $C = \kappa e^{c_\infty +\frac{\xi}{2}R_\X^2}$ under \ref{assump1} and $C = e^{ c_\infty}$ under \ref{assump2} or \ref{assump3}.

Next, recall that
\begin{equation}
\label{eq:strong-concavity-2}
    \frac{\dd^2}{\dd s^2} \Kant(\psi + s v) = \sca{ \frac{\dd^2}{\dd s^2} (\psi + s v)^{c ,\lambda} }{ \mu },
\end{equation}
    and that from the associativity of variances described in \eqref{eq:associativity-var},
    \begin{equation}
    \label{eq:strong-concavity-3}
        - \lambda \frac{\dd^2}{\dd s^2} \Kant(\psi + s v) =  \Var_{\nu[\psi + sv]}(v) - \Var_{ \mu }( \frac{\dd}{\dd s} (\psi + s v)^{c ,\lambda} )  .
    \end{equation}
    Plugging \eqref{eq:strong-concavity-2} and \eqref{eq:strong-concavity-3} into yields \eqref{eq:strong-concavity-4} yields
    \begin{equation*}
       C  \frac{\dd^2}{\dd s^2} \Kant(\psi + s v) \leq - \lambda \frac{\dd^2}{\dd s^2} \Kant(\psi + s v) -  \Var_{\nu[\psi + sv]}(v).
    \end{equation*}
    This gives
    \begin{equation}
      \label{eq:Kantorovich-strong-concavity-prior-to-rescaling}
        \frac{\dd^2}{\dd s^2}
        \Kant(\psi + s v)
        \leq
        -\left(
            C
            +  \lambda
        \right)^{-1}
        \Var_{\nu[\psi + sv]}(v).
    \end{equation}

    We will complete the proof by improving the constant on the right-hand side
    of \eqref{eq:Kantorovich-strong-concavity-prior-to-rescaling} via the
    homogeneity argument described below.
    It will be convenient to introduce the more explicit notation
    that highlights the role played by the regularization parameter $\lambda$
    and the cost function $c$.
    Henceforth, define
    \begin{equation}
      \Kant_{c, \lambda}(\psi) = \sca{\psi^{c, \lambda}}{\mu}
    \end{equation}
    and
    \begin{equation}
      \dd \nu_{c,\lambda}[\psi](y) =
      \dd \nu(y) \int \exp\left(
        \frac{\psi^{c,\lambda}(x) + \psi(y) - c(x,y)}{\lambda}
        \right)
        \dd\mu(x).
    \end{equation}
    Let $\alpha > 0$ be an arbitrary positive constant.
    Observe that for any $\psi$ we have
    \begin{align}
      (\alpha\psi)^{\alpha c, \alpha \lambda}(x)
      &=
      -\alpha\lambda \int
      \exp\left(\frac{\alpha \psi(y) - \alpha c(x,y)}{\alpha\lambda}\right)
      \dd\nu(y)
      \\
      &=
      \alpha
      \cdot
      (\psi)^{c,\lambda}(x).
    \end{align}
    It follows that
    \begin{equation}
      \Kant_{\alpha c, \alpha \lambda}(\alpha \psi)
      = \sca{(\alpha \psi)^{\alpha c, \alpha \lambda}}{\mu}
      = \alpha \Kant_{c, \lambda}(\psi)
    \end{equation}
    Applying \eqref{eq:Kantorovich-strong-concavity-prior-to-rescaling}
    to $\Kant_{\alpha c, \alpha \lambda}$ yields
    \begin{align}
        \frac{\dd^2}{\dd s^2}
        \Kant_{c, \lambda}(\psi + s v)
        &=
        \alpha^{-1}
        \frac{\dd^2}{\dd s^2}
        \Kant_{\alpha c, \alpha \lambda}(\alpha \psi + s \alpha v)
        \\
        &\leq
        -\alpha^{-1}\left(
            C_\alpha
            +  \alpha \lambda
        \right)^{-1}
        \Var_{\nu_{\alpha c, \alpha \lambda}[\alpha \psi + s\alpha v]}(\alpha v),
        \label{eq:Kantorovich-strong-convexity-proof-penultimate-step}
    \end{align}
    where $C_\alpha = \kappa e^{\alpha(c_\infty +\frac{\xi}{2}R_\X^2)}$ if \ref{assump1} holds and
    $C_\alpha = e^{\alpha( c_\infty)}$ if \ref{assump2} or \ref{assump3} hold.
    Finally, observe that for any $\psi$ we have
    \begin{align}
      \dd \nu_{\alpha c, \alpha \lambda}[\alpha \psi](y)
      &=
      \dd \nu(y) \int \exp\left(
        \frac{(\alpha\psi)^{\alpha c,\alpha \lambda}(x) + \alpha \psi(y) -
        \alpha c(x,y)}{\alpha \lambda}
        \right)
        \dd\mu(x)
      \\
      &=
      \dd \nu(y) \int \exp\left(
        \frac{\alpha(\psi)^{c,\lambda}(x) + \alpha \psi(y) -
        \alpha c(x,y)}{\alpha \lambda}
        \right)
        \dd\mu(x)
      \\
      &=
      d\nu_{c, \lambda}[\psi](y).
    \end{align}
    The above identity plugged into
    \eqref{eq:Kantorovich-strong-convexity-proof-penultimate-step} yields
    \begin{align}
        \frac{\dd^2}{\dd s^2}
        \Kant_{c, \lambda}(\psi + s v)
        &\leq
        -\alpha \left(
            C_\alpha
            +  \alpha \lambda
        \right)^{-1}
        \Var_{\nu_{c, \lambda}[\psi + s v]}(v).
        \\
        &=
        -\left(
            \alpha^{-1}C_\alpha
            +  \lambda
        \right)^{-1}
        \Var_{\nu_{c, \lambda}[\psi + s v]}(v).
    \end{align}
    The above inequality holds for any $\alpha > 0$.
    Under assumption \ref{assump1}, taking $\alpha = (c_{\infty} + \frac{\xi}{2}R_{\mathcal{X}}^2)^{-1}$ finishes the proof
    of the first part of the statement in Proposition~\ref{prop:strong-concavity}. Under assumptions \ref{assump2} or \ref{assump3},
    taking $\alpha = (c_{\infty})^{-1}$ finishes the proof of the second part of this statement.

\hfill\qed \newline

\subsection{Proof of Lemma~\ref{lemma:I-concave}}
\label{sec:proof-of-lemma-I-concave}

    The proof of Lemma~\ref{lemma:I-concave} follows from the Prékopa-Leindler inequality \cite{prekopa1, leindler, prekopa2}. We start by recalling the statement of this inequality, found for instance in \cite{cordero2006prekopa-leindler}.
\begin{theorem}[Weighted Prékopa-Leindler inequality (Theorem 1.2 of \cite{cordero2006prekopa-leindler})]
\label{th:prekopa-leindler}
    Let $\xi > 0$ and $\rho$ be a measure on $\Rsp^d$ of the form $\dd \rho = e^{-W}$ where $W$ is a $\xi$-strongly-convex function. Let $ s \in [0,1]$ and let $f,g,h:\Rsp^d \to \Rsp_+$ be such that for all $x,y \in \Rsp^d$,
    $$ h((1- s)x +  s y) \geq e^{-\xi  s(1- s) \nr{x-y}^2/2}f(x)^{1- s} g(y)^{ s}.$$
    Then,
    $$ \int_{\Rsp^d} h \dd \rho \geq \left(\int_{\Rsp^d} f \dd \rho \right)^{1-s} \left(\int_{\Rsp^d} g \dd \rho \right)^{s}. $$
\end{theorem}

 We are able to apply this inequality thanks to the following technical result, which relies on the semi-concavity assumption made on $c$ in \ref{assump}.
    \begin{lemma} \label{lemma:c-lambda-transform-concavity-estimate}
        Let $\psi^0, \psi^1 : \Rsp^d \to \Rsp$, $ s \in (0,1)$, and assume \ref{assump}. Then, denoting $\psi^{ s} = (1 -  s) \psi^0 +  s \psi^1$, it holds for any $u, v \in \X$ that
        \begin{align}
            ( \psi^{ s})^{c, \lambda}((1 -  s)u +  s v) \geq -\frac{\xi}{2}  s (1 -  s) \nr{u-v}^2 + (1- s) (\psi^0)^{c,\lambda}(u) +  s (\psi^1)^{c, \lambda}(v).
        \end{align}
    \end{lemma}
    \begin{proof}[Proof of Lemma~\ref{lemma:c-lambda-transform-concavity-estimate}]
        Let $u, v \in \X$. Assumption \ref{assump} ensures that there exists $\xi \geq 0$ such that for all $y \in \Y$, $x \mapsto c(x, y) - \frac{\xi}{2} \nr{x}^2$ is concave on $\X$. It means in particular, that for any $y \in \Y$,
    \begin{equation*}
        c((1- s)u +  s v, y) - \frac{\xi}{2} \nr{ (1- s)u +  s v }^2 \geq (1- s)( c(u, y) - \frac{\xi}{2} \nr{u}^2) +  s ( c(v, y) - \frac{\xi}{2} \nr{v}^2).
    \end{equation*}
    After rearranging, this ensures that for any $y \in \Y$,
    \begin{equation*}
        c((1- s)u +  s v, y)  \geq (1- s) c(u, y) +  s  c(v, y) - \frac{\xi}{2}  s (1 -  s) \nr{u-v}^2.
    \end{equation*}
    From this inequality together with Hölder's inequality, we deduce:
    \begin{align*}
        ( &(1 -  s) \psi^0 +  s \psi^1)^{c, \lambda}((1 -  s)u +  s v) = - \lambda \log \int e^{ \frac{ (1- s) \psi^0(y) +  s \psi^1(y) - c((1- s)u +  s v, y)}{\lambda} } \dd \nu(y) \\
        &\geq  - \lambda \log \int \left( e^{ \frac{ \psi^0(y) - c(u, y)}{\lambda} } \right)^{1- s} \left( e^{ \frac{ \psi^1(y) - c(v, y)}{\lambda} } \right)^{ s}  e^{\frac{\xi  s(1- s) \nr{u-v}^2}{2\lambda}}\dd \nu(y) \\
        &=- \frac{\xi}{2}  s (1 -  s) \nr{u-v}^2 - \lambda \log \int \left( e^{ \frac{ \psi^0(y) - c(u, y)}{\lambda} } \right)^{1- s} \left( e^{ \frac{ \psi^1(y) - c(v, y)}{\lambda} } \right)^{ s} \dd \nu(y)\\
        &\geq  -\frac{\xi}{2}  s (1 -  s) \nr{u-v}^2  - \lambda \log \left( \int e^{ \frac{ \psi^0(y) - c(u, y)}{\lambda} } \dd \nu(y) \right)^{1- s} \left( \int e^{ \frac{ \psi^1(y) - c(v, y)}{\lambda} } \dd \nu(y) \right)^{ s}   \\
        &= -\frac{\xi}{2}  s (1 -  s) \nr{u-v}^2 + (1- s) (\psi^0)^{c,\lambda}(u) +  s (\psi^1)^{c, \lambda}(v).
    \end{align*} %\qedhere
    \end{proof}

We are now ready to prove Lemma~\ref{lemma:I-concave}.

\begin{proof}[Proof of Lemma~\ref{lemma:I-concave}]
    Let $\psi^0, \psi^1 : \Rsp^d \to \Rsp$, $ s \in (0,1)$, and denote
    \begin{align*}
        h &: u \mapsto \exp( ( (1 -  s) \psi^0 +  s \psi^1)^{c, \lambda}(u) ), \\
        f &: u \mapsto \exp( ( \psi^0 )^{c, \lambda}(u) ), \\
        g &: u \mapsto \exp( ( \psi^1 )^{c, \lambda}(u) ).
    \end{align*}
    Lemma~\ref{lemma:c-lambda-transform-concavity-estimate} guarantees that for any $u, v \in \X$,
    $$ h((1- s)u +  s v) \geq e^{-\xi  s(1- s) \nr{u-v}^2/2} f(u)^{1- s} g(v)^{ s}.$$
    From the weighted Prékopa-Leindler inequality recalled in Theorem~\ref{th:prekopa-leindler}, it thus follows that
    \begin{align*}
        I_\rho ( (1 -  s) \psi^0 +  s \psi^1 ) &= \log \int_\X h(x)\dd \rho(x) \\
        &\geq (1- s) \log \int_\X f(x)\dd \rho(x) +  s \log \int_\X g(x)\dd \rho(x) \\
        &= (1 -  s) I_\rho(\psi^0) +  s I_\rho(\psi^1),
    \end{align*}
    which shows the concavity of $I_\rho$. Note that we used the convexity of $\X$ here.
\end{proof}

\section{Gaussian measures}\label{sec:gaussians}

In this section, we prove a lower-bound on the dual sub-optimality gap along Sinkhorn's iterates in a simple one-dimensional Gaussian setting with linear ground-cost. We also give the value of the limiting contraction rate in the asymptotic regime $\lambda \to 0$.
\begin{theorem}
    \label{th:lower-bound-speed-gaussian}
    Let $\sigma > 0$, $\mu = \mathcal{N}(0, 1)$, $\nu = \mathcal{N}(0, \sigma^2)$, $c(x,y) = -xy$, and $\psi_0 = 0$.
    If 
    $\lambda \leq \sigma/5$, then for any $t \geq 0$ the sub-optimality gap $\delta_{t}$ satisfies 
    $$ \delta_t \geq \frac{2 \sigma^6}{ (\sigma + \lambda)^3} \left(1 - \frac{4 \lambda}{\sigma} - \frac{4 \lambda^2}{\sigma^2}\right)^t(\alpha^*)^2
    \quad\text{where}\quad
     \alpha^* = \frac{\lambda}{4 \sigma^2} - \frac{\sqrt{4 \sigma^2 + \lambda^2}}{4\sigma^2}.
    $$
    Moreover, in the asymptotic regime $\lambda \to 0$, we get the limit
    $$ \lim_{t \to + \infty} \frac{\delta_{t+1}}{\delta_t} = 1 - \frac{4\lambda}{\sigma} + o(\lambda). $$
\end{theorem}
Before proceeding with the proof of the above theorem, we show how it implies Theorem~\ref{thm:lower-bound-gaussian-less-formal} stated in the introduction of this paper.

\begin{proof}[Proof of Theorem~\ref{thm:lower-bound-gaussian-less-formal}]
Observe first that the condition $\lambda \leq \sigma / 5 < \sigma$ entails
$$ \frac{2\sigma^{6}}{(\sigma + \lambda)^3} \geq \frac{\sigma^3}{4}.$$
Then notice that 
$$ (\alpha^*)^2 = \frac{\lambda^2}{16 \sigma^4} + \frac{4 \sigma^2 + \lambda^2}{16 \sigma^4} - \frac{\lambda}{8 \sigma^4}\sqrt{4 \sigma^2 + \lambda^2}.$$
But $\sqrt{4 \sigma^2 + \lambda^2} \leq \sqrt{4 \sigma^2} + \sqrt{\lambda^2} = 2 \sigma + \lambda$, so that
$$ (\alpha^*)^2 \geq \frac{1}{4 \sigma^2} - \frac{\lambda}{4 \sigma^3}.$$
The condition $\lambda \leq \sigma/5$ then ensures that $$(\alpha^*)^2  \geq \frac{1}{5 \sigma^2}.$$ Thus Theorem~\ref{th:lower-bound-speed-gaussian} combined with the two lower-bounds above guarantee that for any $t \geq 0$ we have
$$ \delta_t \geq \frac{2 \sigma^6}{ (\sigma + \lambda)^3} \left(1 - \frac{4 \lambda}{\sigma} - \frac{4 \lambda^2}{\sigma^2}\right)^t(\alpha^*)^2 \geq  \frac{\sigma}{20} \left(1 - \frac{4 \lambda}{\sigma} - \frac{4 \lambda^2}{\sigma^2}\right)^t. $$
Finally, using again that $\lambda \leq \sigma/5$ we get $\frac{4\lambda^2}{\sigma^2} \leq \frac{4 \lambda}{5 \sigma}$, so that 
$$1 - \frac{4 \lambda}{\sigma} - \frac{4 \lambda^2}{\sigma^2} \geq 1 - \frac{24\lambda}{5 \sigma} \geq 1 - \frac{5 \lambda}{\sigma} \geq 0. $$
This yields, for any $t \geq 0$,
$$ \delta_t \geq \frac{\sigma}{20} \left( 1 - \frac{5 \lambda}{\sigma} \right)^t,$$
which proves Theorem~\ref{thm:lower-bound-gaussian-less-formal}. %\hfill\qed
\end{proof}

In order to prove Theorem~\ref{th:lower-bound-speed-gaussian}, we leverage the fact that 
$\mu$ and $\nu$ are both Gaussian to explicitly describe Sinkhorn's iterates. 
Indeed, starting from $\psi_0 = 0$, it can be shown directly that Sinkhorn's iterates are second 
order polynomials of the form 
$$ \psi_t(y) = \alpha_t y^2 + \beta_t \qquad \text{and} \qquad 
\psi_t^{c,\lambda}(x) = \gamma_t x^2 + \omega_t,$$
where $(\alpha_t)_t, (\beta_t)_t, (\gamma_t)_t$ and $(\omega_t)_t$ are sequences of 
real numbers defined recursively as follows:
\begin{align*}
    \begin{cases}
        \alpha_0 = 0, \\
        \beta_0 = 0,
    \end{cases}
    \begin{cases}
        \gamma_t = \frac{\sigma^2}{4 \sigma^2 \alpha_t - 2\lambda}, \\
        \omega_t = -\beta_t + \frac{\lambda}{2} 
        \log\left(1 - \frac{2 \sigma^2}{\lambda} \alpha_t \right),
    \end{cases}
    \begin{cases}
        \alpha_{t+1} = \frac{1}{4\gamma_t - 2 \lambda}, \\
        \beta_{t+1} = -\omega_t + \frac{\lambda}{2} 
        \log\left(1 - \frac{2}{\lambda} \gamma_t \right).
    \end{cases}
\end{align*}
The semi-dual functional 
$E(\psi) = \sca{\psi^{c,\lambda}}{\mu} + \sca{\psi}{\nu}$ satisfies along these 
iterates the equality 
$$ E(\psi_t) = \sigma^2 \alpha_t + \beta_t + \gamma_t + \omega_t. $$
This can be written as a function of $\alpha_t$:
$$E(\psi_t) = \sigma^2 \alpha_t + \frac{\lambda}{2} 
\log\left(1 - \frac{2 \sigma^2}{\lambda} \alpha_t \right) + \frac{\sigma^2}{4 \sigma^2 \alpha_t - 2\lambda}.
$$ 
We thus shift our attention to the sequence $(\alpha_t)_t$. 
The above expressions allow to get 
the recursion
\begin{align}
    \label{eq:def-alpha}
    \begin{cases}
        \alpha_0 = 0, \\
        \alpha_{t+1} = \frac{2 \sigma^2 \alpha_t - \lambda}{2 \sigma^2 + 2 \lambda^2 - 4 \lambda \sigma^2 \alpha_t}.
    \end{cases}
\end{align}
Analysing this recursion allows to get the following convergence result for $(\alpha_t)_t$, 
together with a bound translating the fact that this convergence cannot happen too fast.
\begin{lemma}
    \label{lemma:convergence-alpha}
    The sequence $(\alpha_t)_t$ defined in \eqref{eq:def-alpha} is decreasing and converges to $\alpha^*$ 
    (defined in the statement of Theorem~\ref{th:lower-bound-speed-gaussian}). If $\lambda \leq \sigma/5$, then for any $t \geq 0$ 
    this sequence satifies
    $$ \alpha_t - \alpha^* \geq  \left(1 - \frac{2 \lambda}{\sigma} - \frac{2 \lambda^2}{\sigma^2}\right)^t (-\alpha^*). $$
    Moreover, in the asymptotic regime $\lambda \to 0$, we get the limit
    $$ \lim_{t \to + \infty} \frac{\alpha_{t+1} - \alpha^*}{\alpha_t - \alpha^*} = 1 - \frac{2 \lambda}{\sigma} + o(\lambda). $$
\end{lemma}
\begin{proof}
    We start by noticing that $\alpha_{t+1} = f(\alpha_t)$ where 
    $f : x \mapsto \frac{2 \sigma^2 x - \lambda}{2 \sigma^2 + 2 \lambda^2 - 4 \lambda \sigma^2 x}$. 
    This function satisfies 
    $$ f'(x) = \left( \frac{\sigma^2}{\sigma^2 + \lambda^2 - 2 \lambda \sigma^2 x}\right)^2.$$
    The map $f$ is thus increasing, and by noticing that $\alpha_1 \leq \alpha_0$ we 
    have with an immediate recursion that for all $t \geq 0$, $\alpha_{t+1} \leq \alpha_t$, 
    so that $(\alpha_t)_t$ is a decreasing sequence. Because $\alpha_0 = 0$, this ensures 
    $\alpha_t \leq 0$ for all $t \geq 0$. We can also notice that for $x \leq 0$,
    $$0 \leq f'(x) \leq \left( \frac{\sigma^2}{\sigma^2 + \lambda^2} \right)^2 < 1, $$
    which entails that $f$ is a contraction on $\Rsp_-$. The sequence $(\alpha_t)_t$ 
    thus converges to the unique fixed point of $f$ on $\Rsp_-$, which is 
    $$\alpha^* = \frac{\lambda}{4 \sigma^2} - \frac{\sqrt{4 \sigma^2 + \lambda^2}}{4\sigma^2}.$$
    We differentiate $f$ once more and get
    $$f''(x) = \frac{4 \lambda \sigma^6}{(\sigma^2 + \lambda^2 - 2 \lambda \sigma^2 x)^3}.$$
    Notice that $f'' \geq 0$ on $\Rsp_-$, so that $f$ is convex on this part of $\Rsp$. 
    Because the sequence $(\alpha_t)_t$ lives on $\Rsp_-$, this entails that for any $t\geq 0$ 
    one has 
    $$ f(\alpha_t) \geq f(\alpha^*) + f'(\alpha^*)(\alpha_t - \alpha^*), $$ 
    that is for any $t \geq 0$,
    $$ \alpha_{t+1} - \alpha^* \geq f'(\alpha^*)(\alpha_t - \alpha^*).$$
    An immediate recursion entails that for any $t \geq 0$,
    \begin{equation}
        \label{eq:lower-bound-speed-alpha}
        \alpha_{t} - \alpha^* \geq f'(\alpha^*)^t (\alpha_0 - \alpha^*).
    \end{equation}
    We now seek for a positive lower-bound on $f'(\alpha^*)$. 
    We have $$\alpha^* = \frac{\lambda}{4 \sigma^2} - \frac{\sqrt{4 \sigma^2 + \lambda^2}}{4\sigma^2}.$$
    But $\sqrt{a+b} \leq \sqrt{a} + \sqrt{b}$ for any non-negative reals $a$ and $b$. We thus have the inequality
    $$\alpha^* \geq \frac{\lambda}{4 \sigma^2} - \left(\frac{2 \sigma + \lambda}{4\sigma^2}\right) = - \frac{1}{2 \sigma}.$$
    This yields
    \begin{align*}
        f'(\alpha^*) &= \left( \frac{\sigma^2}{\sigma^2 + \lambda^2 - 2 \lambda \sigma^2 \alpha^*}\right)^2 \\
        &\geq \left( \frac{\sigma^2}{\sigma^2 + \lambda^2 + \lambda \sigma}\right)^2 \\
        &= \left(1 + \frac{\lambda}{\sigma} + \frac{\lambda^2}{\sigma^2} \right)^{-2} \\
        &\geq 1 - \frac{2\lambda}{\sigma} - \frac{2\lambda^2}{\sigma^2}.
    \end{align*}
    Now the assumption $\lambda \leq \frac{\sigma}{5}$ implies that $\lambda < \frac{(\sqrt{3}-1)}{2} \sigma$, which entails in turn that $1 - \frac{2\lambda}{\sigma} - \frac{2\lambda^2}{\sigma^2} > 0$. We thus have for any integer $t \geq 0$
    $$ f'(\alpha^*)^t \geq \left(1 - \frac{2\lambda}{\sigma} - \frac{2\lambda^2}{\sigma^2}\right)^t. $$
    Plugging this last bound into the estimate \eqref{eq:lower-bound-speed-alpha} yields the wanted lower-bound 
    on the speed of convergence of $(\alpha_t)_t$.
    
    We conclude with the limiting contraction rate of the sequence $(\alpha_t)_t$ in the asymptotic regime $\lambda \to 0$. We have, for any $t \geq 0$,
    $$ \frac{\alpha_{t+1} - \alpha_*}{\alpha_t - \alpha^*} = \frac{f(\alpha_{t}) - f(\alpha^*)}{\alpha_t - \alpha^*}.$$
    Since $\alpha_t$ converges to $\alpha^*$ as $t$ goes to infinity, we have the limit
    $$ \lim_{t \to +\infty} \frac{\alpha_{t+1} - \alpha_*}{\alpha_t - \alpha^*} = f'(\alpha^*).$$
    Straightforward computations show that in the asymptotic regime $\lambda \to 0$, it holds
    $$ \alpha^* = -\frac{1}{2 \sigma} + \frac{\lambda}{4 \sigma^2} + o(\lambda).$$
    We thus get the asymptotics
    \begin{align*}
        \lim_{t \to +\infty} \frac{\alpha_{t+1} - \alpha_*}{\alpha_t - \alpha^*} = f'(\alpha^*) 
        &= \left( \frac{\sigma^2}{\sigma^2 + \lambda^2 - 2 \lambda \sigma^2 \alpha^*}\right)^2 \\
        &= \left( 1 + \frac{\lambda}{\sigma} + o(\lambda) \right)^{-2} \\
        &= 1 - \frac{2 \lambda}{\sigma} + o(\lambda),
    \end{align*}
    which is the limit stated in Lemma~\ref{lemma:convergence-alpha}.
\end{proof}

On a side note, notice that Lemma~\ref{lemma:convergence-alpha} allows to recover the formula 
proven in \cite{chen2015optimal, bojilov2016matching, gerolin2019kinetic, janati20-eot-gaussian} for the value of the maximum of the dual functional:
\begin{align*}
    E(\psi_*) &= \frac{\sigma^2}{4 \sigma^2 \alpha^* - 2\lambda} 
+ \frac{\lambda}{2} 
\log\left(1 - \frac{2 \sigma^2}{\lambda} \alpha^* \right) + \sigma^2 \alpha^* \\
    &= \frac{\lambda}{2} \log\left( \frac{1}{2} + \frac{(\lambda^2 + 4\sigma^2)^{1/2}}{2 \lambda} \right) +\frac{\lambda}{2} - \frac{(4 \sigma^2 + \lambda^2)^{1/2}}{2}.
\end{align*}

We now have a bound on the speed of convergence of $(\alpha_t)_t$ towards its limit $\alpha^*$. In 
order to prove Theorem~\ref{th:lower-bound-speed-gaussian}, there remains to translate this bound 
into a bound regarding the convergence of $(E(\psi_t))_t$ towards its limit $E(\psi^*)$, as done bellow.

\begin{proof}[Proof of Theorem \ref{th:lower-bound-speed-gaussian}]
    We assume that $\lambda \leq \frac{\sigma}{5}$ and give ourselves a fixed integer $t \geq 0$. In order to prove the first lower-bound, our objective is 
    to show that
    \begin{equation}
        \label{eq:suboptimality-alphas}
        \delta_t \geq \frac{2\sigma^6}{(\sigma + \lambda)^3} (\alpha_t - \alpha^*)^2. 
    \end{equation}
    Indeed, assuming that \eqref{eq:suboptimality-alphas} holds, Lemma~\ref{lemma:convergence-alpha} ensures
    $$ \delta_t \geq \frac{2\sigma^6}{(\sigma + \lambda)^3} \left(1 - \frac{2 \lambda}{\sigma} - \frac{2 \lambda^2}{\sigma^2}\right)^{2t} (\alpha^*)^2. $$
    But it always holds that $\left(1 - \frac{2 \lambda}{\sigma} - \frac{2 \lambda^2}{\sigma^2}\right)^{2} \geq 1 - \frac{4 \lambda}{\sigma} - \frac{4 \lambda^2}{\sigma^2}$, and from $\lambda \leq \frac{\sigma}{5} \leq \frac{(\sqrt{2}-1)}{2} \sigma$ we have that $1 - \frac{4\lambda}{\sigma} - \frac{4\lambda^2}{\sigma^2} > 0$. This gives
    $$ \delta_t \geq \frac{2\sigma^6}{(\sigma + \lambda)^3} \left(1 - \frac{4 \lambda}{\sigma} - \frac{4 \lambda^2}{\sigma^2}\right)^{t} (\alpha^*)^2, $$
    which is the targeted estimate. We thus now focus on proving \eqref{eq:suboptimality-alphas}.\\
    A repeated use of the fundamental theorem of calculus allows us to write
    \begin{align}
        \delta_t &= E(\psi^*) - E(\psi_t) \\
        &= \int_0^1 \frac{\dd}{\dd s} E(\psi_t + s(\psi^* - \psi_t)) \dd s \\
        &= - \int_0^1 \int_s^1 \frac{\dd^2}{\dd u^2} E(\psi_t + u(\psi^* - \psi_t)) \dd u \dd s + \frac{\dd}{\dd s} E(\psi_t + s(\psi^* - \psi_t))\Bigg|_{s=1}.
    \end{align}
    Because $\psi^*$ is a maximizer of the concave functional $E$, we have the first order condition
    $$ \frac{\dd}{\dd s} E(\psi_t + s(\psi^* - \psi_t))\Bigg|_{s=1} = 0. $$
    Moreover, from the definition of $E$ we have $$\frac{\dd^2}{\dd u^2} E(\psi_t + u(\psi^* - \psi_t)) = \frac{\dd^2}{\dd u^2} \Kant(\psi_t + u(\psi^* - \psi_t)),$$
    where $\Kant$ is the \emph{entropic Kantorovich functional} defined in Section~\ref{sec:background}. The formula given in Lemma~\ref{lemma:K-concave-derivatives} thus ensures
    \begin{align}
        \delta_t = \frac{1}{\lambda} \int_0^1 \int_s^1 \int_\Rsp \Var_{\nu_x[\psi_t + u(\psi^* - \psi_t)]}(\psi_t - \psi_*) \dd \mu(x) \dd u \dd s.
    \end{align}
    Recall now that the iterates $\psi_t$ and $\psi^*$ are second order polynomials satisfying for any $y \in \Rsp$
    $$ \psi_t(y) = \alpha_t y^2 + \beta_t \qquad \text{and} \qquad \psi^*(y) = \alpha^* y^2 + \beta^*. $$
    Hence the sub-optimality reads
    \begin{align}
        \label{eq:expression-suboptimality-gaussian}
        \delta_t = \frac{(\alpha_t - \alpha^*)^2}{\lambda} \int_0^1 \int_s^1 \int_\Rsp \Var_{y \sim \nu_x[\psi_t + u(\psi^* - \psi_t)]}(y^2) \dd \mu(x) \dd u \dd s.
    \end{align}
    In order to ease the notation, introduce for any $u \in [0,1]$ and $x \in \Rsp$ the probability measure $\nu_x^u := \nu_x[\psi_t + u(\psi^* - \psi_t)]$ on $\Rsp$. 
    By definition, the density of 
    this measure satifies for all $y \in \Rsp$
    $$  \nu_x^u(y) = \exp\left( \frac{(\psi_t + u(\psi^* - \psi_t))^{c,\lambda}(x) + (\psi_t + u(\psi^* - \psi_t))(y) + xy}{\lambda} \right) \nu(y). $$
    Introduce now the real numbers 
    \begin{align}\begin{cases}
        \alpha^u := (1-u)\alpha_t + u \alpha^*, \\
         \beta^u := (1-u)\beta_t + u \beta^*,
    \end{cases} \text{and} \quad
    \begin{cases}
        \gamma^u := \frac{\sigma^2}{4 \sigma^2 \alpha^u - 2 \lambda}, \\
        \omega^u :=- \beta^u + \frac{\lambda}{2}\log\left(1 - \frac{2 \sigma^2}{\lambda}\alpha^u\right).
    \end{cases}
\end{align}
They are such that
\begin{align*}
    \forall y \in \Rsp, \quad (\psi_t + u(\psi^* - \psi_t))(y) &= \alpha^u y^2 + \beta^u, \\
    \text{and} \quad \forall x \in \Rsp, \quad 
(\psi_t + u(\psi^* - \psi_t))^{c,\lambda}(x) &= \gamma^u x^2 + \omega^u.
\end{align*}
With these notations we thus have the following expression:
\begin{align*}
    \nu_x^u(y) &= \exp\left( \frac{ \gamma^u x^2 + \omega^u + \alpha^u y^2 + \beta^u + xy}{\lambda} \right) \nu(y) \\
    &=\frac{1}{\sqrt{2\pi}\sigma}\exp\left( \frac{ \gamma^u x^2 + \omega^u + \alpha^u y^2 + \beta^u + xy}{\lambda} - \frac{y^2}{2 \sigma^{2}}\right) \\
    &= \frac{1}{\sqrt{2\pi}\sigma}\exp\left( \frac{ \omega^u + \beta^u }{\lambda}\right) \exp\left( - \frac{1}{2 \hat{\sigma}^{2}}(y + 2\gamma^u x)^2\right), \\
\end{align*}
where $$ \hat{\sigma}^2 = \frac{\lambda \sigma^2}{\lambda - 2 \sigma^2 \alpha^u} = -2 \lambda \gamma^u.$$
Noticing that 
$$  \frac{ \omega^u + \beta^u }{\lambda} = \frac{1}{2}\log\left(1 - \frac{2 \sigma^2}{\lambda}\alpha^u\right) = \frac{1}{2} \log\left( \frac{\sigma^2}{\hat{\sigma}^2}\right), $$
we have the simplification
$$  \nu_x^u(y) = \frac{1}{\sqrt{2\pi}\hat{\sigma}} \exp\left( - \frac{1}{2 \hat{\sigma}^{2}}(y + 2\gamma^u x)^2\right), $$
which corresponds to the density of the Gaussian measure $\mathcal{N}(-2\gamma^u x, \hat{\sigma}^2)$.
We now compute the value of $\Var_{y \sim \nu_x^u}(y^2)$ that appears in \eqref{eq:expression-suboptimality-gaussian}. This variance reads
    $$ \Var_{y \sim \nu_x^u}(y^2) = \int_\Rsp y^4 \dd \nu_x^u(y)  - \left(\int_\Rsp y^2 \dd \nu_x^u(y) \right)^2.$$
Using that $\nu_x^u = \mathcal{N}(-2\gamma^u x, \hat{\sigma}^2)$, we easily get 
$$ \int_\Rsp y^2 \dd \nu_x^u(y) = (2\gamma^u x)^2 + \hat{\sigma}^2 \quad \text{and} 
\quad  \int_\Rsp y^4 \dd \nu_x^u(y) = (2\gamma^u x)^4 + 6 (2\gamma^u x)^2 \hat{\sigma}^2 + 3 \hat{\sigma}^4. $$
Recalling that $ \hat{\sigma}^2 = -2 \lambda \gamma^u,$ this yields
\begin{align*}
    \int_\Rsp y^2 \dd \nu_x^u(y) &= 4(\gamma^u)^2 x^2 - 2\lambda \gamma^u \\
     \text{and} 
\quad  \int_\Rsp y^4 \dd \nu_x^u(y) &= 16(\gamma^u)^4 x^4 -48 (\gamma^u)^3 x^2 + 12 \lambda^2 (\gamma^u)^2. 
\end{align*}
This eventually gives the expression
$$ \Var_{y \sim \nu_x^u}(y^2) =  -32 \lambda (\gamma^u)^3 x^2 + 8 \lambda^2 (\gamma^u)^2. $$
Using that $\mu = \mathcal{N}(0, 1)$, this entails
\begin{equation}
    \label{eq:lowerbound-var-gaussian}
    \int_\Rsp \Var_{y \sim \nu_x^u}(y^2) \dd \mu(x) =  -32 \lambda (\gamma^u)^3  + 8 \lambda^2 (\gamma^u)^2. 
\end{equation}
We now seek a lower bound on this quantity that is independent of $u$. Notice first that $8 \lambda^2 (\gamma^u)^2 > 0$, and that
$$ - (\gamma^u)^3 = \left( \frac{\sigma^2}{2 \lambda - 4 \sigma^2 \alpha^u} \right)^3.$$
But $\alpha^u = (1-u)\alpha_t + u\alpha^* \geq \alpha^*$ from Lemma~\ref{lemma:convergence-alpha}. We have that $\alpha^* \geq - \frac{1}{2 \sigma}$, so that $\alpha^u \geq - \frac{1}{2 \sigma}$ and in turn
$$ - (\gamma^u)^3 \geq \frac{\sigma^6}{8(\sigma + \lambda)^3}.$$
Back to \eqref{eq:lowerbound-var-gaussian}, this yields 
$$ \int_\Rsp \Var_{y \sim \nu_x^u}(y^2) \dd \mu(x) \geq  \frac{4 \lambda \sigma^6}{(\sigma + \lambda)^3}. $$
Injecting this last bound into \eqref{eq:expression-suboptimality-gaussian} finally yields
\begin{align*}
    \delta_t &\geq \frac{4 \sigma^6}{(\sigma + \lambda)^3} \left( \int_0^1 \int_s^1 \dd u \dd s \right) (\alpha_t - \alpha^*)^2 \\
    &= \frac{2 \sigma^6}{(\sigma + \lambda)^3} (\alpha_t - \alpha^*)^2,
\end{align*}
which corresponds to \eqref{eq:suboptimality-alphas}.

We conclude with the proof of the expression of the limit 
$\lim_{t\to +\infty} \frac{\delta_{t+1}}{\delta_t}$ in the asymptotic regime $\lambda \to 0$. From \eqref{eq:expression-suboptimality-gaussian} and \eqref{eq:lowerbound-var-gaussian} we have for any $t \geq 0$ and any $\lambda > 0$ that
$$ \delta_t = \frac{(\alpha_t - \alpha^*)^2}{\lambda} \int_0^1 \int_s^1 \left(  -32 \lambda (\gamma^u)^3  + 8 \lambda^2 (\gamma^u)^2 \right) \dd u \dd s, $$
where $\gamma^u = \frac{\sigma^2}{4 \sigma^2 \alpha^u - 2 \lambda}$ and $\alpha^u = (1-u)\alpha_t + u \alpha^*$. In the asymptotic regime $\lambda \to 0$, it holds
$$ \gamma^u = \frac{1}{4 \alpha^u}\left( 1 + \frac{\lambda}{2 \sigma^2 \alpha^u}\right) + o(\lambda).$$
We thus get in this regime
\begin{align*}
    \delta_t &= \frac{(\alpha_t - \alpha^*)^2}{\lambda} \int_0^1 \int_s^1 \left( - \frac{\lambda}{2 (\alpha^u)^3} \right) \dd u \dd s + o(\lambda) \\
    &= - \frac{1}{4 (\alpha^*)^2} \frac{(\alpha_t - \alpha^*)^2}{\alpha_t} + o(\lambda).
\end{align*}
Hence in the asymptotics $\lambda \to 0$ we obtain the ratio
$$ \frac{\delta_{t+1}}{\delta_t} = \frac{\alpha_t}{\alpha_{t+1}} \left( \frac{\alpha_{t+1} - \alpha^*}{\alpha_t - \alpha^*}\right)^2 + o(\lambda). $$
Taking the limit $t \to + \infty$ of this ratio, Lemma~\ref{lemma:convergence-alpha} ensures that
\begin{align*}
    \lim_{t \to + \infty} \frac{\delta_{t+1}}{\delta_t} &= \frac{\alpha^*}{\alpha^*} \left( 1 - \frac{2 \lambda}{\sigma} \right)^2 + o(\lambda) = 1 - \frac{4 \lambda}{\sigma} + o(\lambda),
\end{align*}
which concludes the proof of Theorem~\ref{th:lower-bound-speed-gaussian}.

\end{proof}

\section{Numerical experiments} \label{sec:numerical-experiments}

This section presents numerical evidence that the core feature of our theoretical findings in Theorems~\ref{thm:quadratic-cost-log-concave} and \ref{thm:quadratic-cost-bounded-density} — the exponential convergence of Sinkhorn's algorithm with a polynomial contraction rate in $\lambda/c_{\infty}$ — also holds in practical settings where both marginals $\mu$ and $\nu$ are discrete probability distributions.

We apply the Sinkhorn's algorithm to compute entropic optimal transport between two grayscale images. The images are represented as probability measures $\mu$ and $\nu$ on a discretized unit square $[0,1]^2$ using an equidistant grid of points; in particular, $c_{\infty} = 1$ and we may focus on the convergence rate as a function of $\lambda$.

\begin{figure}[htbp]
    \centering
    \includegraphics[width=0.85\textwidth]{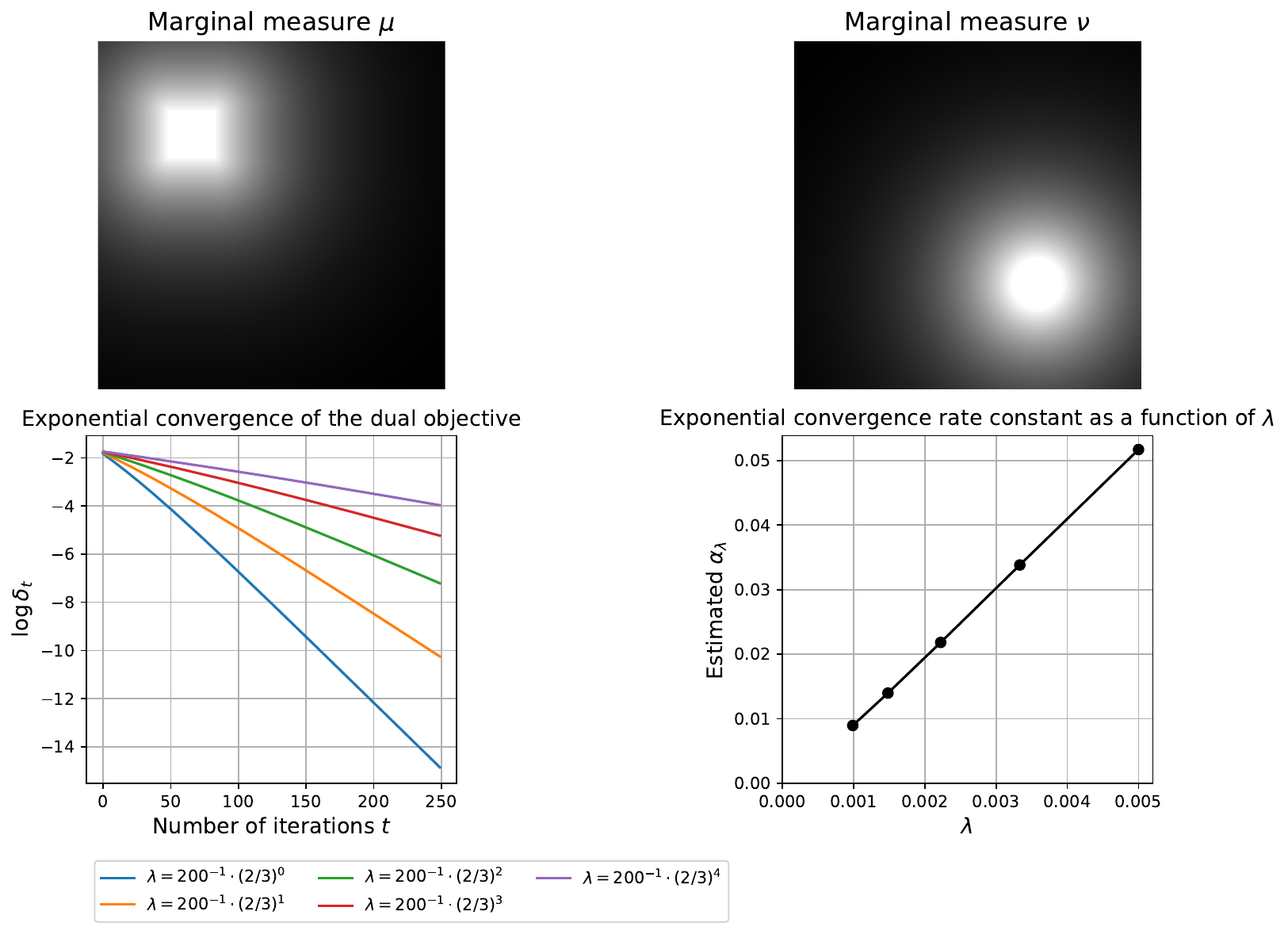}
    \caption{Exponential convergence of Sinkhorn's algorithm with contraction rate $1-\Theta(\lambda/c_{\infty})$ for discretized log-concave marginals.}
    \label{fig:exponential-convergence-log-concave}
\end{figure}

\begin{figure}[htbp]
    \centering
    \includegraphics[width=0.85\textwidth]{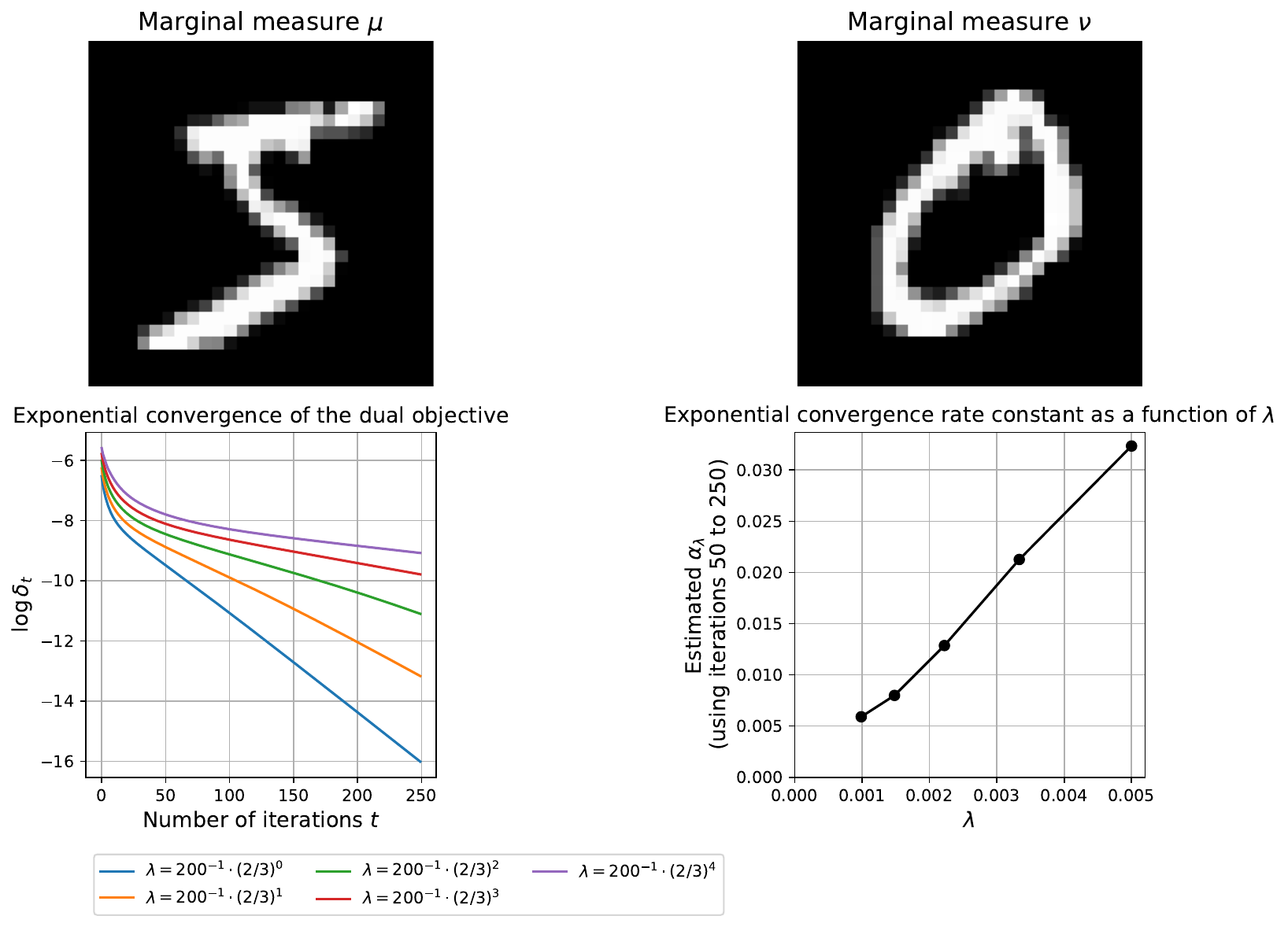}
    \caption{Exponential convergence of Sinkhorn's algorithm for marginal measures corresponding to MNIST images. The contraction rate scales as $1-\Theta(\lambda/c_{\infty})$ for moderately small values of $\lambda$.}
    \label{fig:exponential-convergence-mnist}
\end{figure}

We consider two simulation settings:
\begin{enumerate}
    \item The measures $\mu$ and $\nu$ are two $100 \times 100$ images obtained by discretizing two log-concave distributions. This setting corresponds to that of Theorem~\ref{thm:quadratic-cost-log-concave}, the main difference being that $\mu$ and $\nu$ are discrete measures. See Figure~\ref{fig:exponential-convergence-log-concave} for the simulation results.
    
    \item The measures $\mu$ and $\nu$ are two $28 \times 28$ MNIST images of handwritten digits. Unlike in the first setup, $\mu$ and $\nu$ are not derived from the discretization of underlying log-concave distributions; however, they still exhibit geometric structures typical of real-world images. See Figure~\ref{fig:exponential-convergence-mnist} for the simulation results.
\end{enumerate}

For both simulation settings, and for $\lambda \in \{1/200 \cdot (2/3)^{p} : p \in \{0, 1, \dots, 4\}\}$,  we approximate the optimal dual objective value by running Sinkhorn's algorithm until one-step improvement in the dual objective value reaches machine precision.
Then, for $t \leq 250$, we plot $t \mapsto \log \delta_{t}(\delta)$, where $\delta_{t}(\lambda)$ is the dual suboptimality gap for the dual objective with regularization parameter $\lambda$. Finally, by fitting the least-squares regression line on the observed data-points $\{(t, \log \delta_{t}(\lambda)) : t \leq 250\}$, we find that $\delta_{t}(\lambda) \leq C \cdot (1-\alpha_{\lambda})^{t}$, for $\alpha_{\lambda}$ linear in $\lambda$, in agreement with the main result of our paper proved in Theorem~\ref{thm:quadratic-cost-log-concave}.

The numerical results presented in Figure~\ref{fig:exponential-convergence-log-concave} empirically validate Theorems~\ref{thm:quadratic-cost-log-concave}, \ref{thm:lower-bound-gaussian-less-formal} and the results proved for approximate $(c,\lambda)$-transforms in Appendix~\ref{sec:convergence-rates-approximated-sequence}. On the other hand, when neither of the marginal measures are log-concave nor approximations of log-concave measures, we were only able to prove exponential convergence with contraction rate $1-\Theta(\lambda^2/c_{\infty}^2)$ (see Theorem~\ref{thm:quadratic-cost-bounded-density}). The numerical findings in Figure~\ref{fig:exponential-convergence-mnist} suggest that log-concavity may not be a strict requirement for achieving a contraction rate of order $1 -\Theta(\lambda/c_{\infty})$. This indicates that the main results of this paper might extend to more general settings, which we leave for future work.

\section*{Declarations} No funds, grants, or other support was received.

\bibliographystyle{plain}
\bibliography{ref}

\appendix
\section{Proofs of technical results}
\label{sec:proofs_technical_results}

\subsection{Proof of Lemma~\ref{lemma:K-concave-derivatives}}
\label{sec:proof-of-derivatives-lemma}

First, we state the following lemma on the derivatives of
$(c,\lambda)$-transforms, the proof of which is deferred to the end of the
current section.
\begin{lemma} \label{lemma:derivatives-c-transform}
    Let $\psi, v : \Rsp^d \to \Rsp$. Let $x \in \Rsp^d$. On $\Rsp$, the mapping $\eps \mapsto (\psi + \eps v)^{c, \lambda}(x)$ is $\Class^2$ and admits the following derivatives:
    \begin{align*}
        \frac{\dd}{\dd \eps} (\psi + \eps v)^{c, \lambda}(x) &= -\sca{v}{\nu_x[\psi + \eps v]}, \\
        \frac{\dd^2}{\dd \eps^2}  (\psi + \eps v)^{c, \lambda}(x) &= -\frac{1}{\lambda} \Var_{\nu_x[\psi + \eps v]} (v).
    \end{align*}
\end{lemma}

We are now ready to prove Lemma~\ref{lemma:K-concave-derivatives}.
Let $\psi^0, \psi^1 : \Rsp^d \to \Rsp$ be a pair of potentials. Let $\alpha \in [0, 1]$. We want to show that
    $$ \Kant((1-\alpha) \psi^0 + \alpha \psi^1) \geq (1-\alpha) \Kant(\psi^0) + \alpha \Kant(\psi^1),$$
    with equality if and only is $\psi^0$ and $\psi^1$ differ by an additive constant.
    We have by definition:
    \begin{align*}
        \Kant((1-\alpha) &\psi^0 + \alpha \psi^1) = \int -\lambda \log \int e^{\frac{(1-\alpha) \psi^0(y) + \alpha \psi^1(y) - c(x,y)}{\lambda}} \dd \nu(y) \dd \mu(x) \\
        &= \int -\lambda \log \int \left(e^{\frac{\psi^0(y) - c(x,y)}{\lambda}}\right)^{1-\alpha} \left(e^{\frac{\psi^1(y) - c(x,y)}{\lambda}}\right)^{\alpha}  \dd \nu(y) \dd \mu(x) \\
        &\geq \int -\lambda \log \left (\int e^{\frac{\psi^0(y) - c(x,y)}{\lambda}} \dd \nu(y) \right)^{1-\alpha} \left( \int e^{\frac{\psi^1(y) - c(x,y)}{\lambda}}  \dd \nu(y)\right)^{\alpha}  \dd \mu(x) \\
        &= (1-\alpha) \Kant(\psi^0) + \alpha \Kant(\psi^1),
    \end{align*}
    where the inequality corresponds to Hölder's inequality. There is equality in this inequality if and only if there exists a measurable function $\tau : \Rsp^d \to \Rsp^*_+$ satisfying for $\mu$-a.e. $x$ and for $\nu$-a.e. $y$:
    $$ e^{\frac{\psi^0(y) - c(x,y)}{\lambda}} = \tau(x) e^{\frac{\psi^1(y) - c(x,y)}{\lambda}}.$$
    Simplifying with $e^{\frac{- c(x,y)}{\lambda}}$ on both sides shows that $\tau$ does not depend on $x$, so that equality holds above if and only if there exists a constant $\tau > 0$ such that for $\nu$-a.e. $y$,
    $$ \psi^0(y) = \psi^1(y) + \lambda \log \tau,$$
    that is if and only if $\psi^0$ and $\psi^1$ differ by an additive constant
    $\nu$-a.e. Recalling that $\Kant(\psi) = \sca{(\psi)^{c,\lambda}}{\mu}$,
    the formulas for the first and second derivatives of $\Kant$ follow from
    the formulas for the first and second derivatives of the $(c,
    \lambda)$-transformation given in Lemma~\ref{lemma:derivatives-c-transform}
    and differentiation under the integral sign.
    This completes the proof of Lemma~\ref{lemma:K-concave-derivatives}. \hfill\qed

\begin{proof}[Proof of Lemma~\ref{lemma:derivatives-c-transform}]
    Recall that $$(\psi + \eps v)^{c, \lambda}(x) = -\lambda \log \int \exp \left( \frac{\psi(y) + \eps v(y) - c(x,y)}{\lambda} \right) \dd \nu(y). $$
    This quantity seen as a function of $\eps$ is $\Class^2$ (as a composition of $\Class^2$ functions). It admits the first derivative
    \begin{align*}
        \frac{\dd}{\dd \eps} (\psi + \eps v)^{c, \lambda}(x) &= \frac{ \int -v(y) \exp \left( \frac{\psi(y) + \eps v(y) - c(x,y)}{\lambda} \right) \dd \nu(y) }{\int \exp \left( \frac{\psi(\tilde{y}) + \eps v(\tilde{y}) - c(x,\tilde{y})}{\lambda} \right) \dd \nu(\tilde{y})} \\
        &= -\int v(y) \dd \nu_x[\psi + \eps v](y).
    \end{align*}
    Similarly, we get the second order derivative
    \begin{align*}
        \frac{\dd^2}{\dd \eps^2}&  (\psi + \eps v)^{c, \lambda}(x) \\
        &= \frac{ \int -\frac{v(y)^2}{\lambda} \exp \left( \frac{\psi(y) + \eps v(y) - c(x,y)}{\lambda} \right) \dd \nu(y) }{\int \exp \left( \frac{\psi(\tilde{y}) + \eps v(\tilde{y}) - c(x,\tilde{y})}{\lambda} \right) \dd \nu(\tilde{y})} \\
        &\quad + \frac{ \left(\int v(y) \exp \left( \frac{\psi(y) + \eps v(y) - c(x,y)}{\lambda} \right) \dd \nu(y) \right) \left(\int \frac{v(y)}{\lambda} \exp \left( \frac{\psi(y) + \eps v(y) - c(x,y)}{\lambda} \right) \dd \nu(y) \right) }{ \left( \int \exp \left( \frac{\psi(\tilde{y}) + \eps v(\tilde{y}) - c(x,\tilde{y})}{\lambda} \right) \dd \nu(\tilde{y}) \right)^2 } \\
        &= \frac{-1}{\lambda} \left( \int v(y)^2 \dd \nu_x[\psi + \eps v](y) - \left( \int v(y) \dd \nu_x[\psi + \eps v](y) \right)^2 \right). 
    \end{align*}
    %\hfill\qed
\end{proof}

\subsection{Proof of Corollary \ref{cor:large-t-contraction}}
\label{sec:proof-corollary-large-t-contraction}

Under assumption \ref{assump} and \ref{assump1} or \ref{assump2}, Theorem~\ref{thm:exponential-convergence} and Proposition~\ref{prop:bound-variance-subopt} stated above give constants $\alpha >0$ and $C_1>0$ such that for any $t \geq 0$,
$$ \frac{1}{C_1}  \Var_\nu(\psi^* - \psi_t) \leq \delta_t \leq \left(1 - \alpha^{-1}\right)^t \delta_0. $$
This ensures that the following limit holds:
\begin{equation}
    \label{eq:limit-variance}
    \lim_{t \to +\infty} \Var_\nu(\psi^* - \psi_t) = 0. 
\end{equation}
Now recall that for any $f \in L^2(\nu)$, one has $\Var_\nu(f) = \frac{1}{2} \int \int (f(y) - f(\tilde{y}))^2 \dd \nu(y) \dd \nu(\tilde{y}).$ Thus for any $t \geq 0$,
$$ \Var_\nu(\psi^* - \psi_t) = \frac{1}{2} \nr{\hat{v}_t}^2_{L^2(\nu \otimes \nu)},$$
where 
$$\hat{v}_t : \begin{cases}
    \Y \times \Y &\to \Rsp, \\
    (y ,\tilde{y}) &\mapsto (\psi^* - \psi_t)(y) - (\psi^* - \psi_t)(\tilde{y}).
\end{cases} $$ 
The limit \eqref{eq:limit-variance} can thus be written
\begin{equation}
    \label{eq:limit-variance-rewritten}
    \lim_{t \to + \infty} \nr{\hat{v}_t}_{L^2(\nu \otimes \nu)} = 0.
\end{equation}
We now show that this limit, together with the assumptions made in Corollary~\ref{cor:large-t-contraction}, ensure that
\begin{equation}
\label{eq:uniform-bound-v-hat}
    \lim_{t \to +\infty} \nr{\hat{v}_t}_\infty = \lim_{t \to + \infty} \nr{\psi^* - \psi_t}_{\mathrm{osc}} = 0.
\end{equation}

Indeed, under the assumption that $c$ is Lipschitz continuous, Lemma~\ref{lemma:properties-c-transforms} ensures that the sequence $(\hat{v}_t)_{t \geq 0}$ is uniformly bounded and uniformly continuous. By Arzelà-Ascoli’s theorem, this ensures that this sequence admits a subsequence that uniformly converges to some limit $\hat{v}_\infty \in \Class(\Y \times \Y)$. Because of \eqref{eq:limit-variance-rewritten}, we have that $\nr{\hat{v}_\infty}_{L^2(\nu \otimes \nu)} = 0$, so that $\hat{v}_\infty = 0$ $\nu \otimes \nu$-almost-everywhere. From this we deduce that $\hat{v}_\infty = 0$ on $\Y \times \Y$ (recall that $\Y = \spt(\nu)$ and that $\hat{v}_\infty$ itself is Lipschitz continuous). The limit is thus unique in $\Class(\Y \times \Y)$, therefore the whole sequence $(\hat{v}_t)_{t \geq 0}$ converges uniformly to $0$ in $\Class(\Y \times \Y)$. This corresponds to \eqref{eq:uniform-bound-v-hat}.

Under the assumption that $\nu$ is finitely supported, i.e. that $\Y = \{y_1, \dots, y_N\} \subset \Rsp^d$ for some integer $N \geq 1$, we have
\begin{align*}
    \nr{\hat{v}_t}^2_\infty &= \max_{1 \leq i,j \leq N} \abs{\hat{v}_t(y_i, y_j)}^2 \\ 
    &\leq \frac{1}{\min_{1 \leq i \leq N} \nu(y_i)^2} \sum_{i,j=1}^N \abs{\hat{v}_t(y_i, y_j)}^2 \nu(y_i) \nu(y_j) \\
    &=  \frac{1}{\min_{1 \leq i \leq N} \nu(y_i)^2} \nr{\hat{v}_t}^2_{L^2(\nu \otimes \nu)}. 
\end{align*}
Limit \eqref{eq:uniform-bound-v-hat} then follows from this bound together with \eqref{eq:limit-variance-rewritten}.

Thus in any case, limit \eqref{eq:uniform-bound-v-hat} holds. This ensures that there exists an integer $T \geq 0$, that depends on $\mu, \nu, c$ and $\psi_0$, that is such that for any $t \geq T$,
$$ \nr{\psi^* - \psi_t}_{\mathrm{osc}} \leq \lambda/6. $$
We now assume that $t \geq T$. The above bound, together with the following lemma (the proof of which
can be found in Appendix~\ref{sec:second-order-approximation-lemma-proof}), allow to conclude the proof of Corollary~\ref{cor:large-t-contraction}.

\begin{lemma}
  \label{lemma:semi-dual-second-order-local-approximation}
  Fix any $\psi,v \in L_{1}(\nu)$ such that $\|v\|_{\mathrm{osc}} \leq 1$.
  Then, for any $s \in \Rsp$ such that $|s| \leq \lambda/6$ we have
  \begin{equation}
    E(\psi + sv)
    \leq
    E(\psi) + s\sca{\nu - \nu[\psi]}{v}
  - \frac{s^2}{4}
    \frac{1}{\lambda}
    \Esp_{x \sim \mu} \Var_{\nu_x[\psi]}(v)
  \end{equation}
\end{lemma}
Applying Lemma~\ref{lemma:semi-dual-second-order-local-approximation} with $\psi = \psi^*$, $v = \frac{\psi_t - \psi^*}{\lambda/6}$ and $s = \lambda/6$ yields
$$ E(\psi_{t}) \leq E(\psi^{*}) - \frac{1}{4 \lambda} \Esp_{x \sim \mu} \Var_{\nu_x[\psi^{*}]}(\psi_{t} - \psi^{*}).  $$
Rearranging entails
$$ \Esp_{x \sim \mu} \Var_{\nu_x[\psi^{*}]}(\psi_{t} - \psi^{*}) \leq 4\lambda \delta_{t}. $$
By Lemma~\ref{lemma:K-concave-derivatives} together with Proposition~\ref{prop:strong-concavity}, we have
$$ \Esp_{x \sim \mu} \Var_{\nu_x[\psi^{*}]}(\psi_{t} - \psi^{*}) = -\lambda \frac{\dd^2}{\dd \eps^2} \Kant(\psi^* + \eps v) \Big\vert_{\eps = 0}  \geq\frac{\lambda}{C + \lambda}\Var_{\nu}(v), $$
where $C = e(c_{\infty} + \frac{\xi}{2}R_{\mathcal{X}}^{2})\kappa$ under \ref{assump1} or $C = e c_\infty$ under \ref{assump2}.
The last two bounds give
$$ \Var_{\nu}(v) \leq 4(C + \lambda) \delta_t. $$
Finally, Lemma~\ref{lemma:deriving-contraction-rate} ensures that for $t \geq T$, 
$$ \delta_{t+1} \leq (1-\alpha^{-1})\delta_{t}, $$
with $\alpha = \max\left\{ 64\lambda^{-1}(C+\lambda), \frac{28c_{\infty}}{3}\lambda^{-1} \right\} = 64\lambda^{-1}(C+\lambda)$.

\subsection{Proof of Lemma~\ref{lemma:deriving-contraction-rate}}
\label{sec:proof-of-deriving-contraction-rate-lemma}

Combining Proposition~\ref{prop:one-step-improvement} with the condition
$
  \Var_{\nu}(\psi^{*} - \psi_{t})
  \leq
  C_1 \delta_t + C_2(\delta_t - \delta_{t+1})
$
and applying the inequality $\sqrt{a+b} \leq \sqrt{a} + \sqrt{b}$, valid for
any $a,b > 0$, we have
\begin{align}
    \delta_{t}
    &\leq
    2\sqrt{\lambda^{-1}
      \Var_{\nu}(\psi^{*} - \psi_{t})
      (\delta_{t} - \delta_{t+1})}
    + \frac{14c_{\infty}}{3}\lambda^{-1}(\delta_{t} - \delta_{t+1})
    \\
    &\leq
    2\sqrt{\lambda^{-1}
      (C_1 \delta_t + C_2(\delta_t - \delta_{t+1}))
      (\delta_{t} - \delta_{t+1})}
    + \frac{14c_{\infty}}{3}\lambda^{-1}(\delta_{t} - \delta_{t+1})
    \\
    &=
    2\sqrt{
      \lambda^{-1}C_{1}\delta_{t}(\delta_{t} - \delta_{t+1})
      +
      \lambda^{-1}C_{2}(\delta_{t} - \delta_{t+1})^2
    }
    +
    \frac{14c_{\infty}}{3}\lambda^{-1}(\delta_{t} - \delta_{t+1})
    \\
    &\leq
    2\sqrt{
    \lambda^{-1}C_{1}\delta_{t}(\delta_{t} - \delta_{t+1})}
    +
    \left[
      2\sqrt{\lambda^{-1}C_{2}}
      +
      \frac{14c_{\infty}}{3}\lambda^{-1}
    \right]
    (\delta_{t} - \delta_{t+1})
    \\
    &\leq
    2\max\left\{
      2\sqrt{
      \lambda^{-1}C_{1}\delta_{t}(\delta_{t} - \delta_{t+1})},
      %%%%
      \left[
        2\sqrt{\lambda^{-1}C_{2}}
        +
        \frac{14c_{\infty}}{3}\lambda^{-1}
      \right]
      (\delta_{t} - \delta_{t+1})
    \right\}.
\end{align}
We now split the analysis into two cases, depending on which of the two terms
in the last equation above is larger.

\textbf{Case 1.} We have
\begin{align}
  \delta_{t}
  \leq
  4\sqrt{\lambda^{-1}C_{1}\delta_{t}(\delta_{t} - \delta_{t+1})}.
\end{align}
Squaring both sides and dividing by $\delta_{t} > 0$ we have
\begin{align}
  \delta_{t}
  \leq
  16\lambda^{-1}C_{1}(\delta_{t} - \delta_{t+1}).
\end{align}
Rearranging the above inequality yields
\begin{align}
  \label{eq:contraction-case-1}
  \delta_{t+1} \leq (1-\alpha^{-1})\delta_{t},\quad\text{where}\quad
  \alpha = 16\lambda^{-1}C_{1}.
\end{align}

\textbf{Case 2.}
We have
\begin{align}
  \delta_{t} \leq
    \left[
      4\sqrt{\lambda^{-1}C_{2}}
      +
      \frac{28c_{\infty}}{3}\lambda^{-1}
    \right]
    (\delta_{t} - \delta_{t+1}).
\end{align}
Rearranging the above inequality yields
\begin{align}
  \label{eq:contraction-case-2}
  \delta_{t+1} \leq (1-\alpha^{-1})\delta_{t},\quad\text{where}\quad
  \alpha = 4\sqrt{\lambda^{-1}C_{2}}
      +
      \frac{28c_{\infty}}{3}\lambda^{-1}.
\end{align}

\textbf{Concluding the proof.}
By the two-case analysis outcomes
\eqref{eq:contraction-case-1} and \eqref{eq:contraction-case-2}, we have
\begin{equation}
  \delta_{t+1} \leq (1-\alpha^{-1})\delta_{t},\quad\text{where}\quad
  \alpha =
  \max\left\{
    16\lambda^{-1}C_{1},
    4\sqrt{\lambda^{-1}C_{2}}
      +
      \frac{28c_{\infty}}{3}\lambda^{-1}
    \right\},
\end{equation}
which completes the proof of Lemma~\ref{lemma:deriving-contraction-rate}. \hfill\qed

\subsection{Proof of Lemma~\ref{lemma:variance-comparison-inequality}}
\label{sec:proof-of-variance-comparison-inequality}

Redefining $f$ by $(f - a)/(b-a)$, we may assume without loss of generality
that $a = 0$ and $b = 1$, i.e., $f : \mathcal{X} \to [0,1]$.

The squared Hellinger distance between probability measures $\rho$ and $\pi$ is
defined as
\begin{equation}
  H^{2}(\rho, \pi) = \int_{\mathcal{X}}
  \left(1 - \sqrt{\frac{d\rho}{d\pi}(x)}\right)^2 d\pi(x).
\end{equation}
Using the variational representation of Hellinger distance
\cite[page 149, equation above Example 7.4]{polyanskiy2024information} we have
\begin{equation}
  \label{eq:hellinger-kl-representation}
  \Esp_{\rho}(g)
  \geq
  \frac{1}{2}\Esp_{\pi}(g) - H^{2}(\rho, \pi)
  \quad \forall g : \mathcal{X} \to [0,1].
\end{equation}
Define
\begin{equation}
  m_{\rho} = \Esp_{\rho}(f)
  \quad\text{and}\quad
  h(x) = (f(x) - m_{\rho})^{2}.
\end{equation}
The function $h$ defined above takes values in $[0,1]$; hence,
plugging in $h$ for $g$ in the variational inequality
\eqref{eq:hellinger-kl-representation} yields
\begin{align}
  \Var_{\rho}(f)
  &=
  \Esp_{\rho} (f - m_{\rho})^{2}
  \\
  &\geq
  \frac{1}{2}
  \Esp_{\pi}(f - m_{\rho})^{2}
  - H^{2}(\rho, \pi)
  \\
  &\geq
  \frac{1}{2}\inf_{m \in \mathbb{R}} \Esp_{\pi}(f - m)^2
  - H^{2}(\rho, \pi)
  \\
  &=
  \frac{1}{2}\Var_{\pi}(f)
  - H^{2}(\rho, \pi).
  \label{eq:variance-comparison-hellinger}
\end{align}

Finally, using the inequality $1 - x \leq \log(1/x)$ valid for any $x > 0$
we have
\begin{align}
  H^{2}(\rho, \pi)
  &=
  \int_{\mathcal{X}}
  \left(1 - \sqrt{\frac{d\rho}{d\pi}(x)}\right)^2 d\pi(x)
  \\
  &=
  \int_{\mathcal{X}}
  2\left(1 - \sqrt{\frac{d\rho}{d\pi}(x)}\right) d\pi(x)
  \\
  &\leq
  \int_{\mathcal{X}}
  2\log\left(\sqrt{\frac{d\pi}{d\rho}}\right) d\pi
  \\
  &=
  \int_{\mathcal{X}}
  \log\left(\frac{d\pi}{d\rho}\right) d\pi
  \\
  &=
  \kl{\pi}{\rho}.
\end{align}
Because Hellinger distance is symmetric, it follows that
$$
  H^{2}(\rho, \pi) \leq \min(\kl{\pi}{\rho}, \kl{\rho}{\pi}).
$$
Plugging the above inequality into \eqref{eq:variance-comparison-hellinger}
completes the proof.
\hfill\qed

\subsection{Proof of Proposition~\ref{prop:convexity-sinkhorn-transformation}}
\label{sec:convexity-sinkhorn-transformation}

Let $\psi^0, \psi^1 : \Rsp^d \to \Rsp$ and $\alpha \in [0,1]$. Denote $\psi^\alpha = (1-\alpha)\psi^0 + \alpha \psi^1$. Because assumption \ref{assump} holds, we know from Lemma~\ref{lemma:c-lambda-transform-concavity-estimate} that for any $u, v \in \X$,
$$ ( \psi^\alpha)^{c, \lambda}((1 - \alpha)u + \alpha v) \geq -\frac{\xi}{2} \alpha (1 - \alpha) \nr{u-v}^2 + (1-\alpha) (\psi^0)^{c,\lambda}(u) + \alpha (\psi^1)^{c, \lambda}(v). $$
In turn, the semi-convexity assumption in \ref{assump3} made on $c$ ensures that there exists $\zeta \in \Rsp_+$ such that for all $u, v \in \X$,
$$ c((1-\alpha)u + \alpha v, y) \leq \frac{\zeta}{2} \alpha (1 - \alpha) \nr{u-v}^2 + (1-\alpha) c(u, y) + \alpha c(v, y). $$
%\LC{I think that a factor $\alpha$ is missing in the right-hand side, first term.}
Fix $y \in \Rsp^d$ and denote 
\begin{align*}
        h &: u \mapsto \exp\left( \frac{( \psi^\alpha)^{c, \lambda}(u) - c(u, y)}{\lambda} \right), \\
        f &: u \mapsto \exp\left( \frac{( \psi^0 )^{c, \lambda}(u) - c(u, y)}{\lambda}  \right), \\
        g &: u \mapsto \exp\left( \frac{( \psi^1 )^{c, \lambda}(u) - c(u, y)}{\lambda}  \right).
\end{align*}
The two inequalities above guarantee that for any $u, v \in \X$,
$$ h((1-\alpha)u + \alpha v) \geq e^{-(\frac{\xi+\zeta}{\lambda}) \alpha(1-\alpha) \nr{u-v}^2/2} f(u)^{1-\alpha} g(v)^\alpha.$$
    From the weighted Prékopa-Leindler inequality recalled in Theorem~\ref{th:prekopa-leindler} and the strong convexity assumption made on the function $V$ that is such that $\mu = e^{-V}$ in \ref{assump3}, it thus follows that
    \begin{align*}
        \left( (1-\alpha) \psi_0 + \alpha \psi_1 \right)^{\overline{c, \lambda}}(y) &= -\lambda \log \int_\X h(x)\dd \mu(x) \\
        &\leq -(1-\alpha) \lambda \log \int_\X f(x)\dd \mu(x) - \alpha \lambda \log \int_\X g(x)\dd \mu(x) \\
        &= (1 - \alpha) (\psi^0)^{\overline{c, \lambda}}(y) + \alpha (\psi^1)^{\overline{c, \lambda}}(y),
    \end{align*}
    which shows the convexity of $\psi \mapsto (\psi)^{\overline{c, \lambda}}(y)$ for any $y \in \Rsp^d$. \hfill \qed
\subsection{Proof of
  Lemma~\ref{lemma:semi-dual-second-order-local-approximation}}
\label{sec:second-order-approximation-lemma-proof}
  Notice that:
  \begin{align}
      E(\psi + s v) - E(\psi) - s \sca{v}{\nu - \nu[\psi]} &= E(\psi + s v) - E(\psi) - s \frac{\dd}{\dd \eps} E(\psi + \eps v)\Big|_{\eps = 0} \\
      &= \int_{\eps = 0}^s \frac{\dd}{\dd \eps} E(\psi + \eps v) \dd \eps  - s \frac{\dd}{\dd \eps} E(\psi + \eps v)\Big|_{\eps = 0} \\
      &= \int_{\eps = 0}^s \int_{u=0}^\eps \frac{\dd^2}{\dd u^2} E(\psi + u v) \dd u \dd \eps.
  \end{align}
  Recall that $E(\psi) = \Kant(\psi) + \sca{\psi}{\nu}$, so that $\frac{\dd^2}{\dd u^2} E(\psi + u v) = \frac{\dd^2}{\dd u^2} \Kant (\psi + u v)$. 
  Lemma~\ref{lemma:K-concave-derivatives} thus allows to write
  \begin{equation}
  \label{eq:gap-2nd-derivatives}
      E(\psi + sv) - E(\psi) - s\sca{v}{\nu - \nu[\psi]} 
      = - \frac{1}{\lambda} \int_{\eps = 0}^s \int_{u=0}^\eps \Esp_{x \sim \mu} \Var_{\nu_x[\psi + u v]}(v) \dd u \dd \eps.
  \end{equation} From the definition of $\nu_x[\psi]$ in Section~\ref{sec:background}, we have for any $y \in \Rsp^d$ and $u \in \Rsp_+$ the comparison
  $$ \nu_x[\psi + u v](y) \geq e^{-\frac{2u \nr{v}_\infty}{\lambda}} \nu_x[\psi](y). $$
  The facts that $\nr{v}_\infty \leq 1$ and $s \leq \frac{\lambda}{6}$ entail that for any $y \in \Rsp^d$ and $u \in [0, s]$, 
  $$ \nu_x[\psi + u v](y) \geq e^{-1/3} \nu_x[\psi](y) \geq \frac{1}{2} \nu_x[\psi](y). $$
  This entails in \eqref{eq:gap-2nd-derivatives} the inequality
  \begin{align}
      E(\psi + sv) - E(\psi) - s\sca{v}{\nu - \nu[\psi]} 
      &\leq - \frac{1}{2\lambda} \int_{\eps = 0}^s \int_{u=0}^\eps \Esp_{x \sim \mu} \Var_{\nu_x[\psi]}(v) \dd u \dd \eps \\
      &= - \frac{s^2}{4\lambda} \Esp_{x \sim \mu} \Var_{\nu_x[\psi]}(v),
  \end{align} 
which completes the proof. \hfill\qed
\section{Refining Corollary~\ref{cor:large-t-contraction}}
\label{sec:refinement-cor-value-T}

In this section, we give in specific settings the value of the integer $T$ for which the statement of Corollary~\ref{cor:large-t-contraction} applies.

\begin{proposition}
\label{prop:value-T-c-Lipschitz}
    Assume that \ref{assump} holds, that the cost $c$ is $L$-Lipschitz continuous and
    the target measure $\nu$ is absolutely continuous on a bounded support $\Y$. Assume further that there exists $m_\nu > 0$ such that for any $y \in \Y$ and $r > 0$,
    $$ \nu(B(y, r)) \geq \min(1, m_\nu r^d). $$
    Then, under either \ref{assump1} or \ref{assump2}, the value of $T$ for which Corollary~\ref{cor:large-t-contraction} holds is
    $$ T = \Bigg\lceil\frac{\log\left( (2^4 + \frac{2^{6d+4} L^{2d}}{m_\nu^2}) C_1 \delta_0 \right)}{-\log(1-\alpha^{-1})}\Bigg\rceil.$$ 
    where $C_1$ comes from Proposition~\ref{prop:bound-variance-subopt} and $\alpha$ comes from Theorem~\ref{thm:exponential-convergence}.
\end{proposition}
\begin{remark}
    The regularity condition on the measure $\nu$ is satisfied for instance whenever $\nu$ has a lower-bounded density and $\Y$ is a finite union of convex polytopes, or a finite union of connected domains with smooth boundaries of upper-bounded curvature and all having a lower-bounded Cheeger constant.
\end{remark}

\begin{proposition}
\label{prop:value-T-nu-discrete}
    Assume that \ref{assump} holds and that the target measure $\nu$ is finitely supported on $\Y = \{y_1,\dots,y_N\}$. Then, under either \ref{assump1} or \ref{assump2}, the value of $T$ for which Corollary~\ref{cor:large-t-contraction} holds is
    $$ T = \Bigg\lceil\frac{\log\left(\frac{4 C_1 \delta_0}{\min_{1 \leq i \leq N} \nu(y_i)^2}\right)}{-\log(1-\alpha^{-1})}\Bigg\rceil, $$
    where $C_1$ comes from Proposition~\ref{prop:bound-variance-subopt} and $\alpha$ comes from Theorem~\ref{thm:exponential-convergence}.
\end{proposition}

Propositions \ref{prop:value-T-c-Lipschitz} and \ref{prop:value-T-nu-discrete} 
are both direct consequences of the following lemma, that ensures that the value 
of $T$ in Corollary~\ref{cor:large-t-contraction} can be found provided one 
can control $\nr{\psi_t - \psi^*}_{\mathrm{osc}}$  in terms of 
$\Var_\nu(\psi_t - \psi^*)$.

\begin{lemma}
    \label{lemma:value-T}
    Under assumption \ref{assump} and \ref{assump1} or \ref{assump2}, 
    assume that there exists a constant $C_{2,\infty}$ 
    such that for any $t \geq 0$,
    $$ \nr{\psi_t - \psi^*}^2_{\mathrm{osc}}  \leq C_{2,\infty} 
    \Var_\nu(\psi_t - \psi^*). $$
    Then the value of $T$ such that Corollary \ref{cor:large-t-contraction} holds is 
    $$ T = \Bigg\lceil\frac{\log(2 C_1 C_{2,\infty} \delta_0)}{-\log(1-\alpha^{-1})}\Bigg\rceil, $$
    where $C_1$ is the constant from Proposition~\ref{prop:bound-variance-subopt} and 
    $\alpha$ is the constant from Theorem~\ref{thm:exponential-convergence}.
\end{lemma}
\begin{proof}
    Theorem~\ref{thm:exponential-convergence} together with 
    Proposition~\ref{prop:bound-variance-subopt} ensure that
    $$ \frac{1}{C_1} \Var_\nu(\psi^* - \psi_t) \leq (1 - \alpha^{-1})^t \delta_0. $$
    From our assumption we thus have
    \begin{equation}
        \label{eq:bound-oscillation}
        \nr{\psi_t - \psi^*}^2_{\mathrm{osc}} \leq C_1 C_{2,\infty}  (1 - \alpha^{-1})^t \delta_0. 
    \end{equation}
    Inspecting the proof of Corollary~\ref{cor:large-t-contraction}, the value of $T$ such 
    that this corollary holds is the minimum integer $T \geq 0$ such that for all $t \geq T$,
    $$ \nr{\psi_t - \psi^*}_{\mathrm{osc}} \leq 1. $$
    From \eqref{eq:bound-oscillation}, this is satisfied whenever $t \geq T$ with
    $$ T = \Bigg\lceil\frac{\log(2 C_1 C_{2,\infty} \delta_0)}{-\log(1-\alpha^{-1})}\Bigg\rceil.$$ 
\end{proof}

\begin{proof}[Proof of Proposition~\ref{prop:value-T-c-Lipschitz}]
    Thanks to Lemma~\ref{lemma:value-T}, we only need to find a constant 
    $C_{2,\infty}>0$ that is such that for any $t \geq 0$,
    $$ \nr{\psi_t - \psi^*}^2_{\mathrm{osc}}  \leq C_{2,\infty} 
    \Var_\nu(\psi_t - \psi^*). $$
    As in the proof of Corollary~\ref{cor:large-t-contraction} (Section~\ref{sec:proof-corollary-large-t-contraction}), let us define
    $$\hat{v}_t : \begin{cases}
    \Y \times \Y &\to \Rsp, \\
    (y ,\tilde{y}) &\mapsto (\psi^* - \psi_t)(y) - (\psi^* - \psi_t)(\tilde{y}).
\end{cases} $$ 
This function is such that $\Var_\nu(\psi_t - \psi^*) = \frac{1}{2} \nr{\hat{v}_t}^2_{L^2(\nu \otimes \nu)}$ and $\nr{\psi_t - \psi^*}_{\mathrm{osc}} = \nr{\hat{v}_t}_\infty$.
From Lemma~\ref{lemma:properties-c-transforms}, we know that $\psi_t - \psi^*$ is 
$(2L)$-Lipschitz continous on $\Y$, so that $\hat{v}_t$ is itself $(2L)$-Lipschitz continous 
on $\Y \times \Y$. Let's denote $(y, \tilde{y}) \in \Y \times \Y$ a point that is 
such that
$$ \nr{\hat{v}_t}_\infty = \abs{\hat{v}_t(y, \tilde{y})}.$$
Because $\hat{v}_t$ is $(2L)$-Lipschitz continous, one has for any $(z, \tilde{z}) \in B\left((y, \tilde{y}), \frac{\nr{\hat{v}_t}_\infty}{4L} \right) \cap \Y \times \Y$ that
$$ \abs{\hat{v}_t(z, \tilde{z})} \geq \frac{\nr{\hat{v}_t}_\infty }{2}. $$
Notice then that 
$$ B\left(y, \frac{\nr{\hat{v}_t}_\infty}{8L} \right) \times B\left(\tilde{y}, \frac{\nr{\hat{v}_t}_\infty}{8L} \right) \subset B\left((y, \tilde{y}), \frac{\nr{\hat{v}_t}_\infty}{4L} \right). $$
Therefore we have the following lower-bound:
\begin{align}
    \Var_\nu(\psi_t - \psi^*) &= \frac{1}{2} \int_{\Y \times \Y} \abs{\hat{v}_t(y, \tilde{y})}^2 \dd \nu(y) \dd \nu(\tilde{y}) \\
    &\geq \frac{1}{2} \int_{B\left(y, \frac{\nr{\hat{v}_t}_\infty}{8L} \right) \times B\left(\tilde{y}, \frac{\nr{\hat{v}_t}_\infty}{8L} \right) \cap \Y \times \Y} \frac{\nr{\hat{v}_t}^2_\infty }{4} \dd \nu(y) \dd \nu(\tilde{y}) \\
    &= \frac{1}{8} \nr{\hat{v}_t}^2_\infty \nu\left( B\left(y, \frac{\nr{\hat{v}_t}_\infty}{8L} \right) \right)\nu\left( B\left(\tilde{y}, \frac{\nr{\hat{v}_t}_\infty}{8L} \right)  \right). \label{eq:bound-var-osc-1}
\end{align}
Leveraging the assumptions on $\nu$, we have (because $\Y$ has a Lipschitz boundary) that 
there exists a constant $C_\Y >0$ only depending on $\Y$ such that
\begin{align*}
    \nu\left( B\left(y, \frac{\nr{\hat{v}_t}_\infty}{8L} \right) \right) &\geq \min\left(1, \frac{m_\nu}{2^{3d} L^d} \nr{\hat{v}_t}_\infty^d \right).
\end{align*}
Injecting this last bound into \eqref{eq:bound-var-osc-1} entails that
\begin{equation}
    \Var_\nu(\psi_t - \psi^*) \geq \min\left( \frac{1}{2^3} \nr{\hat{v}_t}_\infty^2, \frac{m_\nu^2}{2^{6d+3} L^{2d}} \nr{\hat{v}_t}_\infty^{2d+2} \right).
\end{equation}
Recalling that $\nr{\hat{v}_t}_\infty = \nr{\psi_t - \psi^*}_{\mathrm{osc}}$, we have either 
$$ \nr{\psi_t - \psi^*}^2_{\mathrm{osc}} \leq 2^3 \Var_\nu(\psi_t - \psi^*), $$
or
$$ \nr{\psi_t - \psi^*}^{2}_{\mathrm{osc}} \leq \frac{1}{c_\infty^{2d}} \nr{\psi_t - \psi^*}^{2d+2}_{\mathrm{osc}} \leq \left( \frac{2^{6d+3} L^{2d}}{m_\nu^2 c_\infty^{2d}} \right)  \Var_\nu(\psi_t - \psi^*).$$
Thus in any case,
$$ \nr{\psi_t - \psi^*}^{2}_{\mathrm{osc}}  \leq \left( 2^3 + \frac{2^{6d+3} L^{2d}}{m_\nu^2 c_\infty^{2d}} \right)  \Var_\nu(\psi_t - \psi^*).$$
Applying Lemma~\ref{lemma:value-T} finally yields the statement of 
Proposition~\ref{prop:value-T-c-Lipschitz}.
\end{proof}

\begin{proof}[Proof of Proposition~\ref{prop:value-T-nu-discrete}]
    Again, let's define
    $$\hat{v}_t : \begin{cases}
        \Y \times \Y &\to \Rsp, \\
        (y ,\tilde{y}) &\mapsto (\psi^* - \psi_t)(y) - (\psi^* - \psi_t)(\tilde{y}).
    \end{cases} $$ 
    This function is such that $\Var_\nu(\psi_t - \psi^*) = \frac{1}{2} \nr{\hat{v}_t}^2_{L^2(\nu \otimes \nu)}$ and $\nr{\psi_t - \psi^*}_{\mathrm{osc}} = \nr{\hat{v}_t}_\infty$.
    As already covered in the proof of Corollary~\ref{cor:large-t-contraction} (Section~\ref{sec:proof-corollary-large-t-contraction}), we have
    $$\nr{\hat{v}_t}^2_\infty \leq  \frac{1}{\min_{1 \leq i \leq N} \nu(y_i)^2} \nr{\hat{v}_t}^2_{L^2(\nu \otimes \nu)}, $$
    that is 
    $$ \nr{\psi_t - \psi^*}^2_{\mathrm{osc}}  \leq \frac{2}{\min_{1 \leq i \leq N} \nu(y_i)^2}
    \Var_\nu(\psi_t - \psi^*). $$ 
Lemma~\ref{lemma:value-T} allows then to conclude.
\end{proof}

\section{Contraction with approximated \texorpdfstring{$(c,\lambda)$}--transforms} 
\label{sec:convergence-rates-approximated-sequence}
In this section, we show that the rates presented in Theorem~\ref{thm:exponential-convergence} extend to the case where the double $(c, \lambda)$-transforms cannot be accurately computed but only approximated with precision $\eps > 0$ in uniform norm. This case covers the practical setting where discrete measures are used in order to approximate the $(c, \lambda)$-transforms appearing in Sinkhorn's algorithm (see Remark~\ref{rk:approx-c-transform}).
\begin{proposition}
\label{prop:convergence-rates-approximated-sequence}
Let $\eps > 0$ be an error parameter. Starting from an arbitrary $\hat{\psi}_0 \in L^\infty(\nu)$, consider a sequence $(\hat{\psi}_t)_{t \geq 0}$ of \emph{approximated} Sinkhorn iterates satisfying for any $t \geq 0$,
$$ \nr{\hat{\psi}_t}_{osc} \leq c_\infty \quad \text{and} \quad \nr{ \hat{\psi}_{t+1} - \hat{\psi}_{t}^{\overline{c, \lambda}} }_\infty \leq \eps. $$
    Then, with the notations of Theorem~\ref{thm:exponential-convergence}, assuming \ref{assump} and either \ref{assump1} or \ref{assump2}, it holds for any $t\geq 0$
    $$ E(\psi^{*}) - E(\hat{\psi}_{t+1}) \leq (1 - \alpha^{-1})(E(\psi^{*}) - E(\hat{\psi}_{t})) + 2\eps. $$
    As a consequence, the \emph{approximated} Sinkhorn iterates satisfy the following contraction up to an error term: for any $t\geq 0$,
    $$ E(\psi^{*}) - E(\hat{\psi}_{t}) \leq (1 - \alpha^{-1})^t(E(\psi^{*}) - E(\hat{\psi}_{0})) + 2\alpha \eps. $$
\end{proposition}
\begin{proof}
    To ease the notations, for any $t \geq 0$ we denote $\hat{\delta}_t = E(\psi^*) - E(\hat{\psi}_t)$ the sub-optimality gap associated with the approximated iterates.
    Let $t \in \mathbb{N}$. Relying on Bernstein’s moment generating function bound for bounded random variables, we can show as in the proof of Proposition~\ref{prop:one-step-improvement} (equation~\eqref{eq:final-corollary-from-bernsteins-bound}) that 
\begin{align}
\label{eq:bound-on-delta-N}
  \hat{\delta}_{t}
  \leq
  \sqrt{2
      \Var_{\nu[\hat{\psi}_{t}]}(\psi^{*} - \hat{\psi}_{t})
      \kl{\nu}{\nu[\hat{\psi}_{t}]}}
  + \frac{2c_\infty}{3}\kl{\nu}{\nu[\hat{\psi}_{t}]},
\end{align}
Lemma~\ref{lemma:variance-comparison-inequality} ensures that
\begin{align*}
    \Var_{\nu[\hat{\psi}_{t}]}(\psi^{*} - \hat{\psi}_{t}) \leq 2 \Var_{\nu}(\psi^{*} - \hat{\psi}_{t}) + 8 c_\infty^2 \kl{\nu}{\nu[\hat{\psi}_{t}]}.
\end{align*}
Using that $\sqrt{a+b} \leq \sqrt{a} + \sqrt{b}$, injecting the above bound in \eqref{eq:bound-on-delta-N} yields
\begin{align}
\label{eq:bernstein-bound-bis}
  \hat{\delta}_{t}
  \leq 2
  \sqrt{
      \Var_{\nu}(\psi^{*} - \hat{\psi}_{t})
      \kl{\nu}{\nu[\hat{\psi}_{t}]}}
  + \frac{14 c_{\infty}}{3}  \kl{\nu}{\nu[\hat{\psi}_{t}]}.
\end{align}
We now bound $\Var_{\nu}(\psi^{*} - \hat{\psi}_{t})$ and $\kl{\nu}{\nu[\hat{\psi}_{t}]}$ in terms of $\hat{\delta}_t$ and $\hat{\delta}_{t+1}$. \\
Let us start with the term  $\Var_{\nu}(\psi^{*} - \hat{\psi}_{t})$. Denote $\hat{v} = \psi^* - \hat{\psi}_t$. Proposition~\ref{prop:strong-concavity} ensures the following lower bound:
    \begin{align}
      \left( C + \lambda \right)
      \hat{\delta}_t
      \geq
      \int_{\eps=0}^1 \int_{s=\eps}^{1} \Var_{\nu[\hat{\psi}_t + s \hat{v}]}(\hat{v}) \dd s \dd \eps,
    \end{align}
    where $C = e(c_{\infty} + \frac{\xi}{2}R_{\mathcal{X}}^2) \kappa$ under \ref{assump1} and $C = ec_\infty$ under \ref{assump2}. It follows by Lemma~\ref{lemma:variance-comparison-inequality} again that
    \begin{align}
      \left( C + \lambda
        \right)
      \hat{\delta}_t
      &\geq
      \int_{\eps=0}^1 \int_{s=\eps}^{1}
      \frac{1}{2}\Var_{\nu}(\hat{v})
      - (2c_{\infty})^2 \kl{\nu}{\nu[\hat{\psi_t} + s\hat{v}]}
      \dd s \dd \eps
      \\
      &=
      \frac{1}{4}\Var_{\nu}(\hat{v})
      -
      4c_{\infty}^2
      \int_{\eps=0}^1 \int_{s=\eps}^{1}
      \kl{\nu}{\nu[\hat{\psi_t} + s\hat{v}]}\dd s \dd \eps.
    \end{align}
    With the exact same computations as in equations \eqref{eq:kl-s-bound}, we have that
    $$ \kl{\nu}{\nu[\hat{\psi_t} + s\hat{v}]} \leq \frac{1-s}{\lambda} \hat{\delta_t}.$$
    The last two equations yield
    \begin{align}
    \label{eq:bound-var-bis}
      \Var_{\nu}(\psi^* - \hat{\psi}^t) &\leq
      \tilde{C} \hat{\delta_{t}}.
 \end{align}
 where $\tilde{C} = 4e(c_{\infty} + \frac{\xi}{2}R_{\mathcal{X}}^2)\kappa + 4 \lambda
        + \frac{8}{3} c_{\infty}^2 \lambda^{-1}$ under~\ref{assump1} and 
 $\tilde{C} = 4ec_{\infty} + 4 \lambda
        + \frac{8}{3} c_{\infty}^2 \lambda^{-1}$ under~\ref{assump2}.\\
Let us now deal with the term $\kl{\nu}{\nu[\hat{\psi}_{t}]}$. Recall that
\begin{align}
    \hat{\delta}_t - \hat{\delta}_{t+1} &= E(\hat{\psi}_{t+1}) - E(\hat{\psi}_{t}) \\
    &= \sca{\hat{\psi}_{t+1}^{c, \lambda}}{\mu} + \sca{\hat{\psi}_{t+1}}{\nu} - \sca{\hat{\psi}_{t}^{c, \lambda}}{\mu} - \sca{\hat{\psi}_t}{\nu} \\
    &= \sca{\hat{\psi}_{t+1}^{c, \lambda} - \hat{\psi}_{t}^{c, \lambda}}{\mu} + \sca{\hat{\psi}_{t+1} - \hat{\psi}_t}{\nu}. \label{eq:upper-bound-improvement-approx}
\end{align}
Remind that by assumption, the approximated Sinkhorn iterates satisfy
$$ \nr{ \hat{\psi}_{t+1} - \hat{\psi}_{t}^{\overline{c, \lambda}} }_\infty \leq \eps. $$
This entails in particular that
$$ \abs{\sca{\hat{\psi}_{t+1} - \hat{\psi}_t}{\nu} - \sca{\hat{\psi}_{t}^{\overline{c, \lambda}} - \hat{\psi}_t}{\nu}} = \abs{\sca{\hat{\psi}_{t+1} - \hat{\psi}_{t}^{\overline{c, \lambda}}}{\nu}} \leq \eps. $$
Similarly, using that the $(c, \lambda)$-transformation is $1$-Lipschitz continuous with respect to the supremum norm (i.e. $\nr{f^{c, \lambda} - g^{c, \lambda}}_\infty \leq \nr{f -g}_\infty$), this assumption entails the bound
$$ \abs{ \sca{\hat{\psi}_{t+1}^{c, \lambda} - \hat{\psi}_{t}^{c, \lambda}}{\mu} - \sca{( \hat{\psi}_{t}^{\overline{c, \lambda}} )^{c, \lambda} - \hat{\psi}_{t}^{c, \lambda}}{\mu} } =  \abs{ \sca{\hat{\psi}_{t+1}^{c, \lambda} - ( \hat{\psi}_{t}^{\overline{c, \lambda}} )^{c, \lambda}}{\mu} } \leq \eps.$$
Plugging the last two bounds into \eqref{eq:upper-bound-improvement-approx}, we obtain
\begin{align}
    \hat{\delta}_t - \hat{\delta}_{t+1} &\geq \sca{( \hat{\psi}_{t}^{\overline{c, \lambda}} )^{c, \lambda} - \hat{\psi}_{t}^{c, \lambda}}{\mu} + \sca{\hat{\psi}_{t}^{\overline{c, \lambda}} - \hat{\psi}_t}{\nu} - 2 \eps \\
    &= F((\hat{\psi}_{t}^{\overline{c, \lambda}} )^{c, \lambda}, \hat{\psi}_{t}^{\overline{c, \lambda}}) - F(\hat{\psi}_{t}^{c, \lambda}, \hat{\psi}_t) - 2\eps.
\end{align}
From the definition of the $(c, \lambda)$-transform, we have
$$ F((\hat{\psi}_{t}^{\overline{c, \lambda}} )^{c, \lambda}, \hat{\psi}_{t}^{\overline{c, \lambda}}) \geq F(\hat{\psi}_{t}^{c, \lambda}, \hat{\psi}_{t}^{\overline{c, \lambda}}). $$
It follows that
\begin{align}
    \hat{\delta}_t - \hat{\delta}_{t+1} &\geq F(\hat{\psi}_{t}^{c, \lambda}, \hat{\psi}_{t}^{\overline{c, \lambda}}) - F(\hat{\psi}_{t}^{c, \lambda}, \hat{\psi}_t) - 2\eps \\
    &= \sca{ \hat{\psi}_{t}^{\overline{c, \lambda}} - \hat{\psi}_t }{\nu} - 2\eps \\
    &= \lambda \kl{\nu}{\nu[\hat{\psi}_{t}]} - 2\eps,
\end{align}
where we used equation~\eqref{eq:KL-nu-nu-psi} in order to get to the last line. We finally have shown the bound
\begin{equation}
    \kl{\nu}{\nu[\hat{\psi}_{t}]} \leq \frac{1}{\lambda}\left(\hat{\delta}_t - \hat{\delta}_{t+1} + 2 \eps \right).
    \label{eq:upper-bound-kl-approx}
\end{equation}
Finally, injecting \eqref{eq:bound-var-bis} and \eqref{eq:upper-bound-kl-approx} into \eqref{eq:bernstein-bound-bis}, we obtain
 \begin{align}
     \hat{\delta}_{t}
  \leq 2
  \sqrt{
     \frac{\tilde{C}}{\lambda}  \hat{\delta}_t \left(\hat{\delta}_t - \hat{\delta}_{t+1} + 2 \eps \right)
      }
  + \frac{14 c_{\infty}}{3\lambda}  \left(\hat{\delta}_t - \hat{\delta}_{t+1} + 2 \eps \right).
 \end{align}
 We may now distinguish between two cases: a first case where 
 $$ 2
  \sqrt{
     \frac{\tilde{C}}{\lambda}  \hat{\delta}_t \left(\hat{\delta}_t - \hat{\delta}_{t+1} + 2 \eps \right)
      } \geq \frac{14 c_{\infty}}{3\lambda}  \left(\hat{\delta}_t - \hat{\delta}_{t+1} + 2 \eps \right),$$
      and a second case with the reverse inequality.
      In the first case, we obtain
$$\hat{\delta}_{t}
  \leq 4
  \sqrt{
     \frac{\tilde{C}}{\lambda}  \hat{\delta}_t \left(\hat{\delta}_t - \hat{\delta}_{t+1} +2 \eps \right) },$$
     which yields after rearranging
     $$ \hat{\delta}_{t+1} \leq \left(1 - \frac{\lambda}{16 \tilde{C}}\right) \hat{\delta}_t + 2 \eps. $$
     Similarly, in the second case described above we get the bound
     $$ \hat{\delta}_{t+1} \leq \left(1 - \frac{3 \lambda}{28 c_\infty}\right) \hat{\delta}_t + 2 \eps. $$
     Recalling that $\tilde{C} = 4e(c_{\infty} + \frac{\xi}{2}R_{\mathcal{X}}^2)\kappa + 4 \lambda
        + \frac{8}{3} c_{\infty}^2 \lambda^{-1}$ under~\ref{assump1} and 
 $\tilde{C} = 4ec_{\infty} + 4 \lambda
        + \frac{8}{3} c_{\infty}^2 \lambda^{-1}$ under~\ref{assump2}, we get in any case that
        $$ 16 \tilde{C} \geq \frac{28c_\infty}{3},$$
        so that we always have
        $$ \hat{\delta}_{t+1} \leq \left(1 - \frac{\lambda}{16 \tilde{C}}\right) \hat{\delta}_t + 2 \eps. $$
        This ends the proof of the first statement of Proposition~\ref{prop:convergence-rates-approximated-sequence}. The second statement of this proposition follows from the first statement and Lemma~\ref{lemma:contraction-of-seq-plus-error} below.      
\end{proof}
The following result is classical and we include its proof for the sake of completeness.
\begin{lemma}
    \label{lemma:contraction-of-seq-plus-error}
        Let $(u_t)_{t\geq 0}$ be a non-negative real sequence. Let $\eta \in (0, 1)$ and $\eps > 0$. If for any $t \in \Nsp$ it holds 
        $$ u_{t+1} \leq (1 - \eta) u_t + \eps,$$
        then for all $t \in \Nsp$,
        $$ u_t \leq (1 - \eta)^t u_0 + \frac{\eps}{\eta}. $$
    \end{lemma}
    \begin{proof}[Proof of Lemma~\ref{lemma:contraction-of-seq-plus-error}]
        Solving the fixed-point equation $x = (1-\eta)x + \eps$ yields $x = \frac{\eps}{\eta}$. Hence for all $t \in \Nsp$ we have
        \begin{align*}
            u_{t+1} - \frac{\eps}{\eta} &\leq (1 - \eta) u_t + \eps - \frac{\eps}{\eta} \\
            &= (1 - \eta) u_t + \eps - (1-\eta) \frac{\eps}{\eta} - \eps \\
            &= (1-\eta) \left(u_t - \frac{\eps}{\eta}\right).
        \end{align*}
        By recursion this yields for all $t \in \Nsp$:
        $$ u_t - \frac{\eps}{\eta} \leq (1-\eta)^t \left( u_0 - \frac{\eps}{\eta} \right)^t \leq  (1-\eta)^t  u_0, $$
        which proves Lemma~\ref{lemma:contraction-of-seq-plus-error}.
    \end{proof}

\end{document}